\documentclass{article}
\usepackage{amsfonts}

\usepackage{amsmath}

\newtheorem{Theorem}{Theorem}[section]
\newtheorem{Definition}[Theorem]{Definition}
\newtheorem{Proposition}[Theorem]{Proposition}
\newtheorem{Lemma}[Theorem]{Lemma}

\newtheorem{Remark}[Theorem]{Remark}

\newtheorem{Hypothesis}[Theorem]{Hypothesis}
\def\qed{\hfill\hbox{\hskip 6pt\vrule width6pt height7pt
depth1pt  \hskip1pt}\bigskip}

\def\E{{{\rm I} \kern -.15em {\rm E}    }}

\newcommand{\reals}{{{\rm I} \kern -.15em {\rm R} }}
\newcommand{\complex}{{{\rm I} \kern -.52em {\rm C} }}
\newcommand{\nat}{{{\rm I} \kern -.15em {\rm N} }}

\def\proof{ {\bf Proof.}  \quad}
\def\endproof{\qquad \vbox to 5.8pt{\offinterlineskip\hrule
        \hbox to 5.8pt{\vrule height 5.1pt\hss\vrule height 5.1pt}\hrule}}

\begin{document}

\date{\today}

\begin{center}
\textbf{\Large Verification Theorems for Stochastic Optimal Control Problems
via a Time Dependent Fukushima - Dirichlet Decomposition}

\vspace{ 5 mm}

{\large Fausto GOZZI}

\vspace{ 3 mm}

{\large Dipartimento di Scienze Economiche e Aziendali}\\[0pt]
{\large Facolt\`a di Economia }\\[0pt]
{\large LUISS - Guido Carli }\\[0pt]
{\large Viale Pola 12, I-00198 Roma, Italy }


\bigskip
{\large Francesco RUSSO}
\vspace{ 3 mm}

{\large Universit\'{e} Paris 13}\\[0pt]
{\large Institut Galil\'{e}e, Math\'{e}matiques}\\[0pt]
{\large 99, av. JB Cl\'{e}ment, F-99430 Villetaneuse, France}
\end{center}

\vspace{ 8 mm}

\textbf{A.M.S. Subject Classification}: Primary: 35R15, 49L10, 60H05, 93B52,
93E20. Secondary: 35J60, 47D03, 47D07, 49L25, 60J35, 93B06, 93B50.

\vspace{ 3 mm}


\vspace{ 3 mm}

\textbf{Abstract.} This paper is devoted to present a method of proving
verification theorems for stochastic optimal control of finite dimensional
diffusion processes without control in the diffusion term. The value
function is assumed to be continuous in time and once differentiable in
the space variable ($C^{0,1}$) instead of once differentiable in time
and twice in space ($C^{1,2}$), like in the classical results.
The results are obtained using a time dependent Fukushima - Dirichlet
decomposition proved in a companion paper by the same authors
 using stochastic calculus via regularization. Applications, examples and
comparison with other similar results are also given.

\section{Introduction\label{INTRO}}

In this paper we want to present a method to get verification theorems for
stochastic optimal control problems of finite dimensional diffusion
processes without control in the diffusion term. The method is based on a
generalized Fukushima - Dirichlet decomposition proved in the companion
paper \cite{FRFGpart1}. Since this Fukushima - Dirichlet decomposition holds
for functions $u:\left[ 0,T\right] \times \mathbb{R}^{n}\rightarrow \mathbb{R%
}$ that are $C^{0}$ in time and $C^{1}$ in space ($C^{0,1}$ in symbols), our
verification theorem has the advantage of requiring less regularity of the
value function $V$ than the classical ones which need $C^{1}$ regularity in
time and $C^{2}$ in space of $V$ ($C^{1,2}$ in symbols), see e.g. \cite[pp.
140, 163, 172]{FS}.

There are also other verification theorems that works in cases when the
value function is nonsmooth: e.g. it is possible to prove a verification
theorem in the case when $V$ is only continuous (see \cite{GSZstochver},
\cite{LYZ}, \cite[Section 5.2]{Zhoubook}) in the framework of viscosity
solutions. However all these results applied to our cases are weaker than
ours, for reasons that are clarified in Section \ref{COMPARISON}.

Since the method is a bit complex and articulated we present first in next
Section \ref{INTROVER} the statement of our verification Theorems \ref
{th:VTintronostrointro} and \ref{th:VTintroinfhor} in a model case with
simplified assumptions (which substantially yield nondegeneracy of the
diffusion coefficient). In the same section we also put Subsection \ref
{INTROOTHER} where we give some comments on the theorem, its applicability
and its relationship with other similar results.

Below we pass to the body of the paper giving first some notations
in Section \ref{NOTATIONS} and then presenting the general
statements (including also possible degeneracy of the diffusion
coefficient) and their proof in Section \ref{PROOFMODEL}. Section
\ref{FEEDBACK} is devoted to necessary conditions and optimal
feedbacks.

Sections \ref{APPL1EXIT} and \ref{APPL2DEG} contain applications of our
technique to more specific classes of problems where other techniques are
more difficult to use. The first is a case of exit time problem where the
HJB equation is nondegenerate but $C^{1,2}$ regularity is not known to hold
due to the lack of regularity of the coefficients; the second is a case
where the HJB equation is degenerate parabolic. Finally in Section \ref
{COMPARISON} we compare our result with other verification techniques.

\section{The statement of the verification theorems in a model case\label%
{INTROVER}}

To clarify our results we describe briefly and informally below the
framework and the statement of the verification theorem that we are going to
prove in a model case. The precise statements and proofs are given in
Section \ref{PROOFMODEL}; then in Sections \ref{APPL1EXIT}, \ref{APPL2DEG}
applications to more specific (and somehow more difficult) classes of
problems are given. We decided this structure since it is difficult to
provide a single general result: we can say that we introduce a technique,
based on the Fukushima-Dirichlet decomposition and on its representation
given in \cite{FRFGpart1}, that can be adapted with some work to different
settings each time with a different adaptation.

First we take a given stochastic basis $\left( \Omega ,\left( \mathcal{F}%
_{s}\right) _{s\geq 0},\mathbb{P}\right) $ that satisfies the so-called
usual conditions, a finite dimensional Hilbert space $A=\mathbb{R}^{n}$ (the
state space), a finite dimensional Hilbert space $E=\mathbb{R}^{m}$ (the
noise space), a set $U\subseteq \mathbb{R}^{k}$ (the control space). We fix
then a terminal time $T\in \left[ 0,+\infty \right] $ (the horizon of the
problem, that can be finite or infinite but is fixed), an initial time and
state $\left( t,x\right) \in \left[ 0,T\right] \times A$ (which will vary as
usual in the dynamic programming approach). The\ state equation is (recall
that $\mathcal{T}_{t}=\left[ t,T\right] \cap \mathbb{R}$)
\begin{eqnarray}
dy(s) &=&\left[ F_{0}\left( s,y\left( s\right) \right) +F_{1}\left(
s,y\left( s\right) ,z\left( s\right) \right) \right] ds+B\left( s,y\left(
s\right) \right) dW(s),\quad s\in \mathcal{T}_{t},  \notag \\
y\left( t\right) &=&x,  \label{eq:stateintro}
\end{eqnarray}
where the following holds (for a matrix $B$ by $\left| \left|
B\right| \right| $ we mean $\sum_{i,j}\left| b_{ij}\right| $ and,
given $E,F$ finite dimensional spaces, by $\mathcal{L}\left(
E,F\right) $ we mean the set of linear operators from $E$ to $F$).

\begin{Hypothesis}
\label{hp:introstateeq}

\begin{enumerate}
\item  $z:\mathcal{T}_{t}\times \Omega \rightarrow U$ (the control process)
is measurable, locally integrable in $t$ for a.e. $\omega \in \Omega $,
adapted to the filtration $\left( \mathcal{F}_{s}\right) _{s\geq 0}$;

\item  $F_{0}:\mathcal{T}_{0}\times \mathbb{R}^{n}\rightarrow \mathbb{R}^{n}$%
, $F_{1}:\mathcal{T}_{0}\times \mathbb{R}^{n}\times U\rightarrow \mathbb{R}%
^{n}$ and $B:\mathcal{T}_{0}\times \mathbb{R}^{n}\rightarrow \mathcal{L}%
\left( E,\mathbb{R}^{n}\right) $ are continuous.

\item  There exists $C>0$ such that $\forall t\in \mathcal{T}_{0}$, $\forall
x_{1},x_{2}\in \mathbb{R}^{n}$,
\begin{eqnarray*}
\left| \left\langle F_{0}\left( t,x_{1}\right) -F_{0}\left(
t,x_{2}\right) ,x_{1}-x_{2}\right\rangle \right| +\left| \left|
B\left( t,x_{1}\right) -B\left( t,x_{2}\right) \right| \right| ^{2}
&\leq &C\left|
x_{1}-x_{2}\right| ^{2}, \\
\left| \left\langle F_{0}\left( t,x\right) ,x\right\rangle \right| +\left|
\left| B\left( t,x\right) \right| \right| ^{2} &\leq &C\left[ 1+\left|
x\right| ^{2}\right]
\end{eqnarray*}
and there exists a constant $K$ such that, $\forall t\in \mathcal{T}_{0}$, $%
\forall x\in \mathbb{R}^{n}$, $\forall z\in U$,
\begin{eqnarray*}
\left| F_{1}\left( t,x_{1},z\right) -F_{1}\left( t,x_{2},z\right) \right|
&\leq &K\left| x_{1}-x_{2}\right| \\
\left| \left\langle F_{1}\left( t,x,z\right) ,x\right\rangle \right| &\leq
&K\left( 1+\left| x\right| +\left| z\right| \right) .
\end{eqnarray*}

\item  For every $x\in \mathbb{R}^{n}$%
\begin{equation*}
\int_{0}^{T}\left[ \left| F_{0}\left( t,x\right) \right| +\left| \left|
B\left( t,x\right) \right| \right| ^{2}\right] dt<+\infty .
\end{equation*}
\end{enumerate}
\end{Hypothesis}

We call $\mathcal{Z}_{ad}\left( t\right) $ the set of admissible control
strategies when the initial time and state is $\left( t,x\right) $ defined
as
\begin{equation*}
\mathcal{Z}_{ad}\left( t\right) =\left\{
\begin{array}{c}
z:\mathcal{T}_{t}\times \Omega \rightarrow U,\quad z\text{ measurable,
adapted to }\left( \mathcal{F}_{s}\right) _{s\geq 0}\text{,} \\
\text{locally integrable in }s\text{ for a.e. }\omega \in \Omega
\end{array}
\right\}
\end{equation*}
and call $y\left( s;t,x,z\right) $ the state process associated with
a given $z\in \mathcal{Z}_{ad}\left( t\right) $; this is the strong
solution of the equation (\ref{eq:stateintro}) which exists and is
unique for any given $ z $ thanks to Theorem 1.2 in
\cite[p.2]{KrylovCetraro}.

Consider now the case $T<+\infty $. We try to minimize the functional
\begin{equation}
J\left( t,x;z\right) =\mathbb{E}\left\{ \int_{t}^{T}l\left( s,y\left(
s;t,x,z\right) ,z\left( s\right) \right) ds+\phi \left( y\left(
T;t,x,z\right) \right) \right\}  \label{eq:cfintro}
\end{equation}
over all control processes $z\in \mathcal{Z}_{ad}\left( t\right) $. On the
coefficients we assume the following.

\begin{Hypothesis}
\label{hp:introcost}$l:\left[ 0,T\right] \times \mathbb{R}^{n}\times
U\rightarrow \mathbb{R}$ and $\phi :\mathbb{R}^{n}\rightarrow \mathbb{R}$
are continuous and such that for each $z\in \mathcal{Z}_{ad}\left( t\right) $
the function $\left( s,\omega \right) \rightarrow l\left( s,y\left( s,\omega
\right) ,z\left( s,\omega \right) \right) $ is lower semiintegrable in $%
\left[ t,T\right] \times \Omega $ (recall that a real valued function $f$ is
lower semiintegrable in a measure space $\left( M,\mu \right) $ if $%
\int_{M}f^{-}d\mu <+\infty $, where $f^{-}$ stands for the negative part of $%
f$ i.e $f^{-}=(-f)\vee 0$).
\end{Hypothesis}

\begin{Remark}
\label{rm:hpcostJbendefinito}The above requirement implies that $J\left(
t,x;z\right) $ is well defined and $>-\infty $ for each $z\in \mathcal{Z}%
_{ad}\left( t\right) $. It is assumed explicitely to cover a more general
class of problems. This is obvious under additional assumptions, e.g. if $l$
and $\phi $ are bounded below. However some interesting problems arising in
economics contain functions $l$ like $-\ln z$ and in this case the lower
semiintegrability of $\left( s,\omega \right) \rightarrow -\ln z\left(
s,\omega \right) $ for each $z\in \mathcal{Z}_{ad}\left( t\right) $ follows
by ad hoc arguments. These are problems with state constraints that are not
explicitely treated in this paper but we want to set our framework so to be
able to treat such class of problems. \textrm{\hfill
\hbox{\hskip 6pt\vrule width6pt
height7pt depth1pt  \hskip1pt}\bigskip }
\end{Remark}

The value function is defined as
\begin{equation}
V\left( t,x\right) =\inf_{z\in \mathcal{Z}_{ad}\left( t\right) }J\left(
t,x;z\right)  \label{eq:defVintro}
\end{equation}
and a control $z^{\ast }\in \mathcal{Z}_{ad}\left( t\right) $ is optimal at $%
\left( t,x\right) $ if $V\left( t,x\right) =J\left( t,x;z^{\ast }\right) $.
The current value Hamiltonian is defined, for $\left( t,x,p,z\right) \in %
\left[ 0,T\right] \times \mathbb{R}^{n}\times \mathbb{R}^{n}\times U$ as
\begin{equation*}
H_{CV}\left( t,x,p;z\right) =\left\langle F_{0}\left( t,x\right)
,p\right\rangle +\left\langle F_{1}\left( t,x,z\right) ,p\right\rangle
+l\left( t,x,z\right)
\end{equation*}
and the (minimum value) Hamiltonian as
\begin{equation*}
H\left( t,x,p\right) =\inf_{z\in U}H_{CV}\left( t,x,p;z\right) .
\end{equation*}
Since the first term of $H_{CV}\left( t,x,p;z\right) $ does not depend on
the control $z$ we will usually define
\begin{equation*}
H_{CV}^{0}\left( t,x,p;z\right) =\left\langle F_{1}\left( t,x,z\right)
,p\right\rangle +l\left( t,x,z\right)
\end{equation*}
and
\begin{equation}
H^{0}\left( t,x,p\right) =\inf_{z\in U}H_{CV}^{0}\left( t,x,p;z\right) ,
\label{eq:H0def}
\end{equation}
so we can write
\begin{equation}
H_{CV}\left( t,x,p;z\right) =:\left\langle F_{0}\left( t,x\right)
,p\right\rangle +H_{CV}^{0}\left( t,x,p;z\right)  \label{eq:HCV0}
\end{equation}
and
\begin{equation}
H\left( t,x,p\right) =:\left\langle F_{0}\left( t,x\right) ,p\right\rangle
+H^{0}\left( t,x,p\right) .  \label{eq:H0}
\end{equation}
When $T<+\infty $ the Hamilton-Jacobi-Bellman (HJB) equation for the value
function is a semilinear parabolic PDE ($\left( t,x\right) \in \left[ 0,T%
\right] \times \mathbb{R}^{n}$)
\begin{equation}
-\partial _{t}v\left( t,x\right) =\frac{1}{2}\mathrm{Tr}\left[ B^{\ast
}\left( t,x\right) \partial _{xx}v(t,x)B\left( t,x\right) \right]
+\left\langle F_{0}\left( t,x\right) ,\partial _{x}v\left( t,x\right)
\right\rangle +H^{0}\left( t,x,\partial _{x}v\left( t,x\right) \right) ,
\label{eq:HJBintropar}
\end{equation}
with the final condition
\begin{equation}
v\left( T,x\right) =\phi \left( x\right) ,\qquad x\in \mathbb{R}^{n}.
\label{eq:HJBintrocondfin}
\end{equation}
To ensure finiteness and continuity of the Hamiltonian we need also to add
the following assumption.

\begin{Hypothesis}
\label{hp:Hamiltonian1}The functions $l$, $F_{1}$ and the set $U$ are such
that the value function is always finite and the Hamiltonian $H^{0}(t,x,p)$
is well defined, finite and continuous for every $\left( t,x,p\right) \in %
\left[ 0,T\right] \times \mathbb{R}^{n}\times \mathbb{R}^{n}$.
\end{Hypothesis}

\begin{Remark}
\label{rm:hpHamiltonian1}\textrm{The above Hypothesis \ref{hp:Hamiltonian1}
is satisfied e.g. when }$l,\phi $\textrm{\ are continuous and }$U$\textrm{\
is compact. Another possibility is to take }$U$\textrm{\ unbounded, }$F_{1}$%
\textrm{\ sublinear, }$l\left( t,x,z\right) =g\left( x\right) +h\left(
z\right) $ \textrm{with }$g$ \textrm{continuous and bounded, }$h$\textrm{\
continuous and such that}
\begin{equation*}
|h(z)|/|z|\longrightarrow +\infty \qquad as\;\;|z|\rightarrow +\infty ;
\end{equation*}
\textrm{this case was studied e.g in \cite{CDPSiamHJ2,CPDE,EFJDE} even in
the infinite dimensional case (see on this the Sections \ref{APPL1EXIT}, \ref
{APPL2DEG}).\hfill
\hbox{\hskip 6pt\vrule
width6pt height7pt depth1pt \hskip1pt}\bigskip }
\end{Remark}

The statement of the classical verification theorem for this model problem
when $T<+\infty $ is the following, see Definition \ref{df:solHJBstrict} for
the definition of strict solution of equation (\ref{eq:HJBintropar})-(\ref
{eq:HJBintrocondfin}).

\begin{Theorem}
\label{th:VTintroclassical}Assume that Hypotheses
\ref{hp:introstateeq}, \ref {hp:introcost}, \ref{hp:Hamiltonian1}
hold true. Let $v \in C^{1,2}\left( \left[ 0,T\right] \times
\mathbb{R}^{n}\right)$  be a polynomially growing strict solution of
the HJB equation (\ref{eq:HJBintropar})-(\ref {eq:HJBintrocondfin})
on $\left[ 0,T\right] \times \mathbb{R}^{n}$. Then the two following
properties hold true.

\begin{enumerate}
\item[(i)]  $v\leq V$ on $\left[ 0,T\right] \times \mathbb{R}^{n}$;

\item[(ii)]  fix $\left( t,x\right) \in \left[ 0,T\right] \times \mathbb{R}%
^{n}$; if $z\in \mathcal{Z}_{ad}\left( t\right) $ is such that, calling $%
y\left( s\right) =y\left( s;t,x,z\right),$%
\begin{equation*}
H^{0}\left( s,y\left( s\right) ,\partial _{x}v\left( s,y\left( s\right)
\right) \right) =H_{CV}^{0}\left( s,y\left( s\right) ,\partial _{x}v\left(
s,y\left( s\right) \right) ,z\left( s\right) \right) ,\quad \mathbb{P-}a.s.,
\end{equation*}
for a.e. $s\in \left[ t,T\right]$, then $z$ is optimal at $\left( t,x\right)
$ and $v\left( t,x\right) =V\left( t,x\right) $.
\end{enumerate}
\end{Theorem}

This theorem states a sufficient optimality condition and its proof is based
on It\^{o}'s formula, see e.g. \cite[p.268]{Zhoubook}. In the classical
context also necessary conditions and existence of optimal feedback can be
proved, see again \cite[p.268]{Zhoubook}.

The statement of our Verification Theorem in the model problem described
above is very similar. We give it in the case when the following
nondegeneracy hypothesis hold
\begin{equation}
B^{-1}\left( t,x\right) F_{1}\left( t,x,z\right) \text{ is bounded on }\left[
0,T\right] \times \mathbb{R}^{n}\times U,  \label{eq:hpGirsanovintro}
\end{equation}
leaving the discussion for the general case to Section \ref{PROOFMODEL}
(strong solutions of the HJB equation (\ref{eq:HJBintropar})-(\ref
{eq:HJBintrocondfin}) are defined in Definition \ref{df:solHJBstrong}).

\begin{Theorem}
\label{th:VTintronostrointro}Assume that Hypotheses
\ref{hp:introstateeq}, \ref{hp:introcost}, \ref{hp:Hamiltonian1} and
(\ref{eq:hpGirsanovintro}) hold true. Let $v \in C^{0,1}\left(
\left[ 0,T\right] \times \mathbb{R}^{n}\right)$  be a polynomially
growing strong solution of the HJB equation
(\ref{eq:HJBintropar})-(\ref {eq:HJBintrocondfin}) on $\left[
0,T\right] \times \mathbb{R}^{n}$. Then (i) and (ii) of Theorem
\ref{th:VTintroclassical} above hold.
\end{Theorem}

In fact we will go further:

\begin{itemize}
\item  proving a result in the case when (\ref{eq:hpGirsanovintro}) does not
hold;

\item  showing that in our setting we can obtain a necessary condition and
the existence of optimal feedback controls on the same line of the classical
results, see Section \ref{FEEDBACK}.
\end{itemize}

Moreover in a paper in preparation we will show that also some cases where
the drift is a distribution can be treated and also cases where the solution
of HJB enjoys weaker regularity.

When $T=+\infty $, $F_{0}$, $F_{1}$ and $B$ do not depend on time, we
consider the problem of minimizing
\begin{equation*}
J\left( t,x;z\right) =\mathbb{E}\left\{ \int_{t}^{+\infty }e^{-\lambda
t}l_{1}\left( y\left( s;t,x,z\right) ,z\left( s\right) \right) ds\right\} ,
\end{equation*}
where $l\left( t,x,z\right) =e^{-\lambda t}l_{1}\left( x,z\right) $, with $%
\lambda >0$, $l_{1}$ continuous and bounded. In this case the value function
$V$ is defined as in (\ref{eq:defVintro}) but its dependence on $t$ becomes
trivial as $V\left( t,x\right) =e^{\lambda t}V\left( 0,x\right) $ for every $%
t$ and $x$. Then we set $t=0$ and call $V_{0}\left( x\right) =V\left(
0,x\right) $. The HJB equation for $V_{0}\left( x\right) $ becomes an
elliptic PDE ($x\in \mathbb{R}^{n}$)
\begin{equation}
\lambda v\left( x\right) =\frac{1}{2}\mathrm{Tr}\left[ B^{\ast }\left(
x\right) \partial _{xx}v(x)B\left( x\right) \right] +\left\langle
F_{0}\left( t,x\right) ,\partial _{x}v\left( t,x\right) \right\rangle
+H^{1}\left( x,\partial _{x}v\left( t,x\right) \right) ,
\label{eq:HJBintroell}
\end{equation}

where
\begin{equation}
H^{1}\left( x,p\right) =\inf_{z\in U}\left\{ \left\langle F_{1}\left(
x,z\right) ,p\right\rangle +l_{1}\left( x,z\right) \right\} =:\inf_{z\in
U}H_{CV}^{1}\left( x,p;z\right) .  \label{eq:H1def}
\end{equation}

In this elliptic case the statement of our verification theorem is the
following (for the definition of strong solution of the HJB equation (\ref
{eq:HJBintroell}) see Definition \ref{df:solHJBstrongell}).

\begin{Theorem}
\label{th:VTintroinfhor}Assume that Hypotheses \ref{hp:introstateeq} and (%
\ref{eq:hpGirsanovintro}) hold. Assume also that $l_{1}$ is continuous and
bounded and that Hypothesis \ref{hp:Hamiltonian1} hold true with $H^{1}$ in
place of $H^{0}$. If $v$\ is a bounded strong solution of the HJB equation (%
\ref{eq:HJBintroell}) and $v\in C^{1}\left( \mathbb{R}^{n}\right) $\ then $%
v\leq V_{0}$. Moreover fix $x\in \mathbb{R}^{n}$; if $z\in \mathcal{Z}%
_{ad}\left( 0\right) $ is such that, calling $y\left( s\right) =y\left(
s;0,x,z\right) $,
\begin{equation*}
H^{1}\left( y\left( s\right) ,\partial _{x}v\left( s,y\left( s\right)
\right) \right) =H_{CV}^{1}\left( y\left( s\right) ,\partial _{x}v\left(
s,y\left( s\right) \right) ,z\left( s\right) \right) ,\quad \text{for a.e. }%
s\geq 0,\quad \mathbb{P-}a.s.,
\end{equation*}
then $z$ is optimal at $\left( 0,x\right) $ and $v\left( x\right)
=V_{0}\left( x\right) $.\smallskip
\end{Theorem}

\subsection{Other Verification Theorems for stochastic optimal control
problems\label{INTROOTHER}}

We recall that in the literature on stochastic optimal control other
verification theorems for not $C^{1,2}$ functions were proved. Roughly
speaking we can say that each of them (and the technique of proof, too) is
strictly connected with the concept of weak solution of the HJB equation
that is considered. We have in fact the following.

\begin{itemize}
\item  \emph{(Strong solutions).} In \cite{CPDE,EFJDE} (see also \cite
{SCSIAM1,DPD1,GM}) in a setting similar to ours but in infinite dimension, a
verification theorem is proved assuming that there exists a strong solution $%
v$ of HJB which is $C^{1}$ in space. In fact the method outlined here goes
in the same direction, but improves and generalizes in the finite
dimensional case the ideas contained there. Such improvement is not a
straightforward one, as we need to use completely different tools to get our
results and we get more powerful theorems. We refer to Sections \ref
{APPL1EXIT} and \ref{APPL2DEG} for examples and to Section \ref{COMPARISON}
for explanations.

\item  \emph{(Viscosity solutions).} In \cite{LYZ,Zhoubook,GSZstochver} a
verification technique for viscosity solutions is introduced and studied.
Such technique adapts to the case when $v$ is only continuous and so is very
general and applicable to cases when the control enters in the diffusion
coefficients $B$. However in the case of our interest, i.e. when $v$ is $%
C^{1}$ in space, such technique gives weaker results as it requires more
assumptions on the coefficients $F_{0}$, $F_{1}$, $B$ and on the candidate
optimal strategy; see Section \ref{COMPARISON} for explanations.

\item  \emph{(Mild solutions).} In \cite[Theorem 7.2]{FTHJB} a verification
theorem is given when the HJB equation admits mild solutions which are $%
C^{1} $ in space (see Definition \ref{df:mildHJBappl1} or Remark \ref
{rm:betterthanalternative}) and when the solution can be represented using
the solution of a suitable backward SDE. The results available with this
technique are limited to infinite dimensional cases where mild solutions
exist and Girsanov theorem can be applied in a suitable way. Such
restrictions prevent the use of this technique (at least to the present
state of the art) e.g. in the cases described in Sections \ref{APPL1EXIT}
and \ref{APPL2DEG}, see Section \ref{COMPARISON} for explanations.
\end{itemize}

Moreover we want to stress the fact that the Fukushima-Dirichlet
decomposition we got in \cite{FRFGpart1} is in fact stronger than what we
need to prove the verification theorems above. In fact for such purpose it
would be enough to prove a Dynkin-type formula (roughly speaking an It\^{o}
formula after expectation) which is in general easier. We decided to state
and prove such Fukushima-Dirichlet decomposition since it can help to deal
with stochastic control problems where the criterion to minimize does not
contain expectation (e.g. pathwise optimality and optimality in
probability). In such context the HJB equation becomes a stochastic PDE and
to get a verification theorem, Dynkin-type formulae cannot be used. See e.g.
\cite{DPDMT,LS1,LS2,Peng} on this subject.

To sum up we think that the interest of our verification results is the fact
that

\begin{itemize}
\item  we obtain them using only the fact that the solution $v$ of HJB
belongs to $C^{0,1}$ and not necessarily to $C^{1,2}$ obtaining better
results than other techniques.

\item  they can be applied also to problems with pathwise optimality and
optimality in probability.
\end{itemize}

\begin{Remark}
As a first step of our work, we consider here the case of a stochastic
optimal control problem in finite dimension with no state constraints. We
are aware of the fact that in many cases, the HJB equation associated with
the control problem admits a $C^{1,2}$ solution, so that our method does not
give a real advantage in this case. However

\begin{itemize}
\item  in some cases, mainly degenerate ones (see Sections \ref{APPL1EXIT}
and \ref{APPL2DEG}) it is known that the solution of HJB is $C^{0}$ (or $%
C^{\alpha }$ with H\"{o}lder exponent $\alpha >0$) in time and $C^{1}$ in
space but $C^{1,2}$ regularity is not known at this stage. In particular we
show an explicit example where the value function is $C^{1}$ but cannot be $%
C^{2}$ in space;

\item  we think that such technique could be extended also in cases where
HJB equation is intended in a generalized sense; in a paper in preparation,
we investigate a verification theorem related to an SDE whose drift is the
derivative of a continuous function therefore a Schwartz distribution, see
for instance \cite{frw}. That state equation models a stochastic particle
moving in a random (irregular) medium;

\item  we consider problems without state constraints; again we think that
some state constraints problems can be treated with our approach but we do
not handle them here: we provide an example with exit time which usually
presents similar difficulties;

\item  we think that our method can be extended to the infinite dimensional
case, where $C^{1,2}$ regularity results for the solution of HJB
equation are much less known while $C^{1}$ regularity can be found
in a variety of cases (see e.g.
\cite{CDPSiamHJ2,CPDE,EFJDE,SCSIAM1}). We do not perform this
extension here leaving it for further work.\hfill\qed
\end{itemize}
\end{Remark}

\begin{Remark}
Similar ideas as in this work have been expressed with a different
formalism in \cite{chitmania}. There the authors have implemented a
generalized It\^{o} formula in the Krylov spirit for proving an
existence and uniqueness result for a generalized solution of
Bellman equation for controlled processes, in a non-degeneracy
situation.\hfill\qed
\end{Remark}

\section{Notations\label{NOTATIONS}}

Throughout this paper we will denote by $\left( \Omega ,\mathcal{F},\mathbb{P%
}\right) $ a given stochastic basis, where $\mathcal{F}$ stands for a given
filtration $\left( \mathcal{F}_{s}\right) _{s\geq 0}$ satisfying the usual
conditions. Given a finite dimensional real Hilbert space $E$, $W$ will
denote a cylindrical Brownian motion with values in $E$ and adapted to $%
\left( \mathcal{F}_{s}\right) _{s\geq 0}$. Given $0\leq t\leq T\leq +\infty $
and setting $\mathcal{T}_{t}=\left[ t,T\right] \cap \mathbb{R}$ the symbol $%
\mathcal{C}_{\mathcal{F}}\left( \mathcal{T}_{t}\times \Omega ;E\right) $,
will denote the space of all continuous processes adapted to the filtration $%
\mathcal{F}$ with values in $E$. This is a Fr\'{e}chet space if endowed with
the topology of the uniform convergence in probability ($u.c.p.$ from now
on). To be more precise this means that, given a sequence $\left(
X^{n}\right) \subseteq \mathcal{C}_{\mathcal{F}}\left( \mathcal{T}_{t}\times
\Omega ;E\right) $ and $X\in \mathcal{C}_{\mathcal{F}}\left( \mathcal{T}%
_{t}\times \Omega ;E\right) $ we have
\begin{equation*}
X^{n}\rightarrow X
\end{equation*}
if and only if for every $\varepsilon >0$, $t_{1}\in \mathcal{T}_{t}$%
\begin{equation*}
\lim_{n\rightarrow +\infty }\sup_{s\in \left[ t,t_{1}\right] }\mathbb{P}%
\left( \left| X_{s}^{n}-X_{s}\right| _{E}>\varepsilon \right) =0.
\end{equation*}

Given a random time $\tau \geq t$ and a process $(X_{s})_{s\in \mathcal{T}%
_{t}},$ we denote by $X^{\tau }$ the stopped process defined by $X_{s}^{\tau
}=X_{s\wedge \tau }.$ The space of all processes in $\left[ t,T\right] $,
adapted to $\mathcal{F}$ and square integrable with values in $E$ is denoted
by $L_{\mathcal{F}}^{2}\left( t,T;E\right) $. $S^{n}$ will denote the space
of all symmetric matrices of dimension $n.$

Let $k\in \mathbb{N}$. As usual $C^{k}\left( \mathbb{R}^{n}\right) $ is the
space of all functions $:\mathbb{R}^{n}\rightarrow \mathbb{R}$ that are
continuous together with their derivatives up to the order $k$. This is a Fr%
\'{e}chet space equipped with the seminorms
\begin{equation}
\sup_{x\in K}\left| u\left( x\right) \right| _{\mathbb{R}}+\sup_{x\in
K}\left| \partial _{x}u\left( x\right) \right| _{\mathbb{R}^{n}}+\sup_{x\in
K}\left| \partial _{xx}u\left( x\right) \right| _{\mathbb{R}^{n\times n}}+...
\label{eq:normK}
\end{equation}
for every compact set $K\subset \subset \mathbb{R}^{n}$. This space will be
denoted simply by $C^{k}$ when no confusion may arise. If $K$ ia compact
subset of $\mathbb{R}^{n}$ then $C^{k}\left( K\right) $ is a Banach space
with the norm (\ref{eq:normK}). The symbol $C_{b}^{k}\left( \mathbb{R}%
^{n}\right) $ will denote the Banach space of all functions from $\mathbb{R}%
^{n}$ to $\mathbb{R}$ that are continuous and bounded together with their
derivatives up to the order $k$. This space is endowed with the usual $\sup $
norm. Passing to parabolic spaces we denote by $C^{0}\left( \mathcal{T}%
_{t}\times \mathbb{R}^{n}\right) $ the space of all functions
\begin{equation*}
u:\mathcal{T}_{t}\times \mathbb{R}^{n}\rightarrow \mathbb{R},\quad \quad
\left( s,x\right) \mapsto u\left( s,x\right)
\end{equation*}
that are continuous. This space is a Fr\'{e}chet space equipped with the
seminorms
\begin{equation*}
\sup_{\left( s,x\right) \in \left[ t,t_{1}\right] \times K}\left| u\left(
s,x\right) \right| _{\mathbb{R}}
\end{equation*}
for every $t_{1}>0$ and every compact set $K\subset \subset \mathbb{R}^{n}$%
). Moreover we will denote by $C^{1,2}\left( \mathcal{T}_{t}\times \mathbb{R}%
^{n}\right) $ (respectively $C^{0,1}\left( \mathcal{T}_{t}\times \mathbb{R}%
^{n}\right) $), the space of all functions
\begin{equation*}
u:\mathcal{T}_{t}\times \mathbb{R}^{n}\rightarrow \mathbb{R},\quad \quad
\left( s,x\right) \mapsto u\left( s,x\right)
\end{equation*}
that are continuous together with their derivatives $\partial _{t}u$, $%
\partial _{x}u$, $\partial _{xx}u$ (respectively $\partial _{x}u$). This
space is a Fr\'{e}chet space equipped with the seminorms
\begin{eqnarray*}
&&\sup_{\left( s,x\right) \in \left[ t,t_{1}\right] \times K}\left| u\left(
s,x\right) \right| _{\mathbb{R}}+\sup_{\left( s,x\right) \in \left[ t,t_{1}%
\right] \times K}\left| \partial _{s}u\left( s,x\right) \right| _{\mathbb{R}%
^{n}} \\
&&+\sup_{\left( s,x\right) \in \left[ t,t_{1}\right] \times K}\left|
\partial _{x}u\left( s,x\right) \right| _{\mathbb{R}^{n}}+\sup_{\left(
s,x\right) \in \left[ t,t_{1}\right] \times K}\left| \partial _{xx}u\left(
s,x\right) \right| _{\mathbb{R}^{n\times n}}
\end{eqnarray*}
(respectively
\begin{equation*}
\sup_{\left( s,x\right) \in \left[ t,t_{1}\right] \times K}\left| u\left(
s,x\right) \right| _{\mathbb{R}}+\sup_{\left( s,x\right) \in \left[ t,t_{1}%
\right] \times K}\left| \partial _{x}u\left( s,x\right) \right| _{\mathbb{R}%
^{n}}\text{)}
\end{equation*}
for every $t_{1}>0$ and every compact set $K\subset \subset \mathbb{R}^{n}$.
This space will be denoted simply by $C^{1,2}$ (respectively $C^{0,1}$) when
no confusion may arise.

Similarly, for $\alpha ,\beta \in \left[ 0,1\right] $ one defines $C^{\alpha
,1+\beta }\left( \mathcal{T}_{t}\times \mathbb{R}^{n}\right) $ (or simply $%
C^{\alpha ,1+\beta }$) as the subspace of $C^{0,1}\left( \mathcal{T}%
_{t}\times \mathbb{R}^{n}\right) $ of functions $u:\mathcal{T}_{t}\times
\mathbb{R}^{n}\mapsto \mathbb{R}$ such that are $u\left( \cdot ,x\right) $
is $\alpha -$H\"{o}lder continuous and $\partial _{x}u\left( s,\cdot \right)
$ is $\beta -$H\"{o}lder continuous (with the agreement that $0$-H\"{o}lder
continuity means just continuity). If such properties hold just locally then
such space is denoted by $C_{loc}^{\alpha ,1+\beta }\left( \mathcal{T}%
_{t}\times \mathbb{R}^{n}\right) $. Similarly, given an open subset of $%
\mathcal{O}$ of $\mathbb{R}^{n}$ one can define $C^{\alpha ,1+\beta }\left(
\mathcal{T}_{t}\times \mathcal{O}\right) $ and $C_{loc}^{\alpha ,1+\beta
}\left( \mathcal{T}_{t}\times \mathcal{O}\right) $. If $K$ is a compact
subset of $\mathcal{T}_{t}\times \mathbb{R}^{n}$ we define $C^{k}\left(
K\right) $ or $C^{\alpha ,1+\beta }\left( K\right) $ as above. Similarly to $%
C_{b}^{k}\left( \mathbb{R}^{n}\right) $ we define the Banach spaces $%
C_{b}^{0}\left( \mathcal{T}_{t}\times \mathbb{R}^{n}\right) $ $%
C_{b}^{1,2}\left( \mathcal{T}_{t}\times \mathbb{R}^{n}\right) $, $%
C_{b}^{\alpha ,1+\beta }\left( \mathcal{T}_{t}\times \mathbb{R}^{n}\right) $%
, $C_{b}^{0,1}\left( \mathcal{T}_{t}\times \mathbb{R}^{n}\right) $.

\section{Proof of the verification theorems in the model case\label%
{PROOFMODEL}}

In this section we give the precise statement and the proof of the
verification Theorems \ref{th:VTintronostrointro} and \ref{th:VTintroinfhor}
for the model problem described in Section \ref{INTROVER}. Then in Sections
\ref{APPL1EXIT} and \ref{APPL2DEG}\ we will consider two families of
problems to which our technique apply.

Consider first the parabolic case fixing the horizon $T<+\infty $. Under
Hypotheses \ref{hp:introstateeq}, \ref{hp:introcost} and \ref
{hp:Hamiltonian1} we seek to minimize the cost $J\left( t,x;z\right) $ given
in (\ref{eq:cfintro}) over all admissible $z\in \mathcal{Z}_{ad}\left(
t\right) $ where the state equation is given by (\ref{eq:stateintro}).

We define the operator
\begin{equation*}
\mathcal{L}_{0}:D\left( \mathcal{L}_{0}\right) \subseteq C^{0}\left( \left[
0,T\right] \times \mathbb{R}^{n}\right) \longrightarrow C^{0}\left( \left[
0,T\right] \times \mathbb{R}^{n}\right) ,\qquad D\left( \mathcal{L}%
_{0}\right) =C^{1,2}\left( \left[ 0,T\right] \times \mathbb{R}^{n}\right) ,
\end{equation*}
\begin{equation*}
\mathcal{L}_{0}v\left( t,x\right) =\partial _{t}v\left( t,x\right)
+\left\langle F_{0}\left( t,x\right) ,\partial _{x}v\left( t,x\right)
\right\rangle +\frac{1}{2}\text{Tr }\left[ B^{\ast }\left( t,x\right)
\partial _{xx}v\left( t,x\right) B\left( t,x\right) \right] .
\end{equation*}
The HJB equation associated with the problem (\ref{eq:stateintro}) - (\ref
{eq:cfintro}) can then be written as
\begin{equation}
\mathcal{L}_{0}v\left( t,x\right) +H^{0}\left( t,x,\partial _{x}v\left(
t,x\right) \right) =0,\quad v\left( T,x\right) =\phi \left( x\right) ,
\label{eq:HJBristretta}
\end{equation}
where $H^{0}$ is given in (\ref{eq:H0def}).

Now we want to apply the representation result proved in \cite{FRFGpart1}
Section 4, which we recall below for the reader's convenience. First we
consider the following Cauchy problem for $h\in C^{0}\left( \left[ 0,T\right]
\times \mathbb{R}^{n}\right) $.
\begin{equation}
\mathcal{L}_{0}u\left( s,x\right) =h\left( s,x\right) ,\qquad u\left(
T,x\right) =\phi \left( x\right) ,  \label{eq:CPlinear}
\end{equation}
with the following definitions of solution.

\begin{Definition}
\label{df:solstrict}We say that $u\in C^{0}\left( \left[ 0,T\right] \times
\mathbb{R}^{n}\right) $ is a strict solution to the backward Cauchy problem (%
\ref{eq:CPlinear}) if $u\in D\left( \mathcal{L}_{0}\right) $ and (\ref
{eq:CPlinear}) holds.
\end{Definition}

\begin{Definition}
\label{df:solstrong}We say that $u\in C^{0}\left( \left[ 0,T\right] \times
\mathbb{R}^{n}\right) $ is a strong solution to the backward Cauchy problem (%
\ref{eq:CPlinear}) if there exists a sequence $\left( u_{n}\right) \subset
D\left( \mathcal{L}_{0}\right) $ and two sequences $\left( \phi _{n}\right)
\subseteq C^{0}\left( \mathbb{R}^{n}\right) $, $\left( h_{n}\right)
\subseteq C^{0}\left( \left[ 0,T\right] \times \mathbb{R}^{n}\right) $, such
that

\begin{enumerate}
\item  For every $n\in \mathbb{N}$ $u_{n}$ is a strict solution of the
problem
\begin{equation*}
\mathcal{L}_{0}u_{n}\left( t,x\right) =h_{n}\left( t,x\right) ,\qquad
u_{n}\left( T,x\right) =\phi _{n}\left( x\right) .
\end{equation*}

\item  The following limits hold
\begin{equation*}
\begin{array}{l}
u_{n}\longrightarrow u\text{ in }C^{0}\left( \left[ 0,T\right] \times
\mathbb{R}^{n}\right) , \\
h_{n}\longrightarrow h\text{ in }C^{0}\left( \left[ 0,T\right] \times
\mathbb{R}^{n}\right) , \\
\phi _{n}\longrightarrow \phi \text{ in }C^{0}\left( \mathbb{R}^{n}\right) .
\end{array}
\end{equation*}
\end{enumerate}
\end{Definition}

The representation result is the following (Corollaries 4.6 and 4.8 of \cite
{FRFGpart1}).

\begin{Theorem}
\label{th:Itostrong}Let
\begin{equation*}
b_{1}:\left[ 0,T\right] \times \mathbb{R}^{n}\times \Omega \rightarrow
\mathbb{R}^{n},
\end{equation*}
be a continuous progressively measurable field (continuous in $(s,x)$) and
\begin{equation*}
b:\left[ 0,T\right] \times \mathbb{R}^{n}\rightarrow \mathbb{R}^{n},\qquad
\sigma :\left[ 0,T\right] \times \mathbb{R}^{n}\rightarrow \mathcal{L}\left(
\mathbb{R}^{m},\mathbb{R}^{n}\right) ,
\end{equation*}
be continuous functions. Let $u\in C^{0,1}\left( \left[ 0,T\right] \times
\mathbb{R}^{n}\right) $ be a strong solution of the Cauchy problem (\ref
{eq:CPlinear}).

Fix $t\in \left[ 0,T\right] $, $x\in \mathbb{R}^{n}$ and let $(S_{s})$ be a
solution to the SDE
\begin{equation*}
dS_{s}=b_{1}\left( s,S_{s}\right) ds+\sigma \left( s,S_{s}\right)
dW_{s},\qquad S_{t}=x.
\end{equation*}
Then
\begin{eqnarray*}
u\left( s,S_{s}\right) &=&u\left( t,S_{t}\right) +\int_{t}^{s}h\left(
r,S_{r}\right) dr+\int_{t}^{s}\left\langle \partial _{x}u\left(
r,S_{r}\right) ,b_{1}\left( r,S_{r}\right) -b\left( r,S_{r}\right)
\right\rangle dr \\
&&+\int_{0}^{s}\partial _{x}u\left( r,S_{r}\right) \sigma \left(
r,S_{r}\right) dW_{r}.
\end{eqnarray*}
provided that: either we can choose the approximating sequence $\left(
u_{n}\right) $ of Definition \ref{df:solstrong} so that for every $0\leq
t\leq s\leq T$%
\begin{equation}
\lim_{n\rightarrow +\infty }\int_{t}^{s}\left\langle \partial
_{x}u_{n}\left( r,S_{r}\right) -\partial _{x}u\left( r,S_{r}\right)
,b_{1}\left( r,S_{r}\right) -b\left( r,S_{r}\right) \right\rangle
dr=0,\qquad u.c.p.,  \label{eq:hpnewconvderucp}
\end{equation}
or the function
\begin{equation*}
\left( t,x,\omega \right) \rightarrow \sigma ^{-1}\left( t,x\right) \left[
b_{1}\left( r,x,\omega \right) -b\left( r,x\right) \right] ,
\end{equation*}
(where $\sigma ^{-1}$ stands for the pseudo-inverse of $\sigma $), is well
defined and bounded on $\left[ 0,T\right] \times \mathbb{R}^{n}\times \Omega
$.
\end{Theorem}

\begin{Remark}
\label{rm:suffperconvder}If
\begin{equation*}
\lim_{n\rightarrow +\infty }\partial _{x}u_{n}=\partial _{x}u,\qquad \text{%
in }C^{0}\left( \left[ 0,T\right] \times \mathbb{R}^{n}\right),
\end{equation*}
then Assumption (\ref{eq:hpnewconvderucp}) is verified. This means that the
result of Theorem \ref{th:Itostrong} above applies if we know that $u$ is a
strong solution in a more restrictive sense, i.e. substituting the point 2
of Definition \ref{df:solstrong} with
\begin{equation*}
\begin{array}{l}
u_{n}\longrightarrow u\text{ in }C^{0}\left( \left[ 0,T\right] \times
\mathbb{R}^{n}\right), \\
\partial _{x}u_{n}\longrightarrow \partial _{x}u\text{ in }C^{0}\left( \left[
0,T\right] \times \mathbb{R}^{n}\right), \\
h_{n}\longrightarrow h\text{ in }C^{0}\left( \left[ 0,T\right] \times
\mathbb{R}^{n}\right), \\
\phi _{n}\longrightarrow \phi \text{ in }C^{0}\left( \mathbb{R}^{n}\right) .
\end{array}
\end{equation*}
This is a particular case of our setting and it is the one used e.g in \cite
{CPDE,EFJDE} to get the verification result. We can say that in these papers
a result like Theorem \ref{th:Itostrong} is proved under the assumption that
$u$ is a strong solution in this more restrictive sense. It is worth to note
that in such simplified setting the proof of Theorem \ref{th:Itostrong}
follows simply by using standard convergence arguments. In particular there
one does not need to use the Fukushima-Dirichlet decomposition presented in
Section 3 of \cite{FRFGpart1}. So, from the methodological point of view
there is a serious difference with the result of Theorem \ref{th:Itostrong},
see Section \ref{COMPARISON} for comments.\hfill
\hbox{\hskip 6pt\vrule width6pt height7pt depth1pt
\hskip1pt}\bigskip
\end{Remark}

To apply Theorem \ref{th:Itostrong}\ to our case we need first to adapt the
notion of strong solution and to rewrite the main assumptions. Let us give
the following definitions.

\begin{Definition}
\label{df:solHJBstrict}We say that $v\in C^{0}\left( \left[ 0,T\right]
\times \mathbb{R}^{n}\right) $ is a strict solution of the HJB equation (\ref
{eq:HJBristretta}) if $v\in D\left( \mathcal{L}_{0}\right) $ and (\ref
{eq:HJBristretta}) holds on $\left[ 0,T\right] \times \mathbb{R}^{n}$.
\end{Definition}

\begin{Definition}
\label{df:solHJBstrong}A function $v\in C^{0,1}\left( \left[ 0,T\right]
\times \mathbb{R}^{n}\right) $ is a strong solution of the HJB equation (\ref
{eq:HJBristretta}) if, setting $h\left( t,x\right) =-H^{0}\left(
t,x,\partial _{x}v\left( t,x\right) \right) $, $v$ is a strong solution of
the backward Cauchy problem
\begin{equation*}
\mathcal{L}_{0}v\left( t,x\right) =h\left( t,x\right) ,\qquad v\left(
T,x\right) =\phi \left( x\right) ,
\end{equation*}
in the sense of Definition \ref{df:solstrong}.
\end{Definition}

We now present our verification theorem in the extended version, as
announced in Section \ref{INTROVER}. We need the following assumption.

\begin{Hypothesis}
\label{hp:exsolHJBvera}There exists a function $v\in C^{0,1}\left( \left[ 0,T%
\right] \times \mathbb{R}^{n}\right) $ which is a strong solution of
the HJB equation (\ref {eq:HJBristretta}) in the sense of Definition
\ref{df:solHJBstrong} and is polynomially growing with its space
derivative in the variable $x$. Moreover:

\begin{itemize}
\item[(i)]  either we can choose the approximating sequence $\left(
v_{n}\right) $ of Definition \ref{df:solHJBstrong} so that for every $0\leq
t\leq s\leq T$ and for every admissible control $z\in \mathcal{Z}_{ad}\left(
t\right) $%
\begin{equation*}
\lim_{n\rightarrow +\infty }\int_{t}^{s}\left\langle \partial
_{x}v_{n}\left( r,y\left( r\right) \right) -\partial _{x}v\left( r,y\left(
r\right) \right) ,F_{1}\left( r,y\left( r\right) ,z\left( r\right) \right)
\right\rangle dr=0,\qquad u.c.p.,
\end{equation*}

\item[(ii)]  or the function
\begin{equation*}
\left( t,x,z\right) \rightarrow B^{-1}\left( t,x\right) F_{1}\left(
t,x,z\right) ,
\end{equation*}
where $B^{-1}$ stands for the pseudo-inverse of $B$, is well defined and
bounded on $\left[ 0,T\right] \times \mathbb{R}^{n}\times U$.
\end{itemize}
\end{Hypothesis}

\begin{Remark}
\label{rm:singat0proofmodel}All the results below still hold true with
suitable modifications if we assume that:
\begin{itemize}
\item  the strong solution $v$ belongs to $ C^{0}\left( \left[ 0,T\right]
\times \mathbb{R}^{n}\right) \cap C^{0,1}\left( \left[ \varepsilon ,T\right]
\times \mathbb{R}^{n}\right) $ for every small $\varepsilon >0$;

\item  for some $\beta \in \left( 0,1\right) $ the map $\left( t,x\right)
\rightarrow t^{\beta }\partial _{x}v\left( t,x\right) $ belongs to $%
C^{0}\left( \left[ 0,T\right] \times \mathbb{R}^{n}\right) $.
\end{itemize}

The proof of Theorem \ref{th:VTintronostro} in this case is a
straightforward generalization of the one presented here: we do not
give it here to avoid technicalities since here we deal with a model
problem. In Section \ref{APPL1EXIT} we will treat a case with such
difficulty. See also Remark 4.10 of \cite{FRFGpart1} on
this.\hfill\qed
\end{Remark}

We give now here the precise statement of Theorem \ref{th:VTintronostrointro}%
.

\begin{Theorem}
\label{th:VTintronostro}Assume that Hypotheses \ref{hp:introstateeq}, \ref
{hp:introcost}, \ref{hp:Hamiltonian1} and \ref{hp:exsolHJBvera} hold. Let $%
H_{CV}^{0}$, $H^{0}$ be as in (\ref{eq:H0def})-(\ref{eq:HCV0}). Let $v\in
C^{0,1}\left( \left[ 0,T\right] \times \mathbb{R}^{n}\right) $ be a strong
solution of (\ref{eq:HJBristretta}) and fix $\left( t,x\right) \in \left[ 0,T%
\right] \times \mathbb{R}^{n}$ which is polynomially growing with its space
derivative in the variable $x$. Then

\begin{enumerate}
\item[(i)]  $v\leq V$ on $[0,T]\times \mathbb{R}^{n}$.

\item[(ii)]  If $z$ is an admissible control at $\left( t,x\right) $ that
satisfies (setting $y\left( s\right) =y\left( s;t,x,z\right) $)
\begin{equation}
H^{0}\left( s,y\left( s\right) ,\partial _{x}v\left( s,y\left( s\right)
\right) \right) =H_{CV}^{0}\left( s,y\left( s\right) ,\partial _{x}v\left(
s,y\left( s\right) \right) ;z\left( s\right) \right) ,  \label{eq:condsuff}
\end{equation}
for $a.e.$ $\,s\in \left[ t,T\right] $, $\mathbb{P-}$almost surely, then $z$
is optimal at $\left( t,x\right) $ and $v\left( t,x\right) =V\left(
t,x\right) .$
\end{enumerate}
\end{Theorem}

The proof of this theorem follows by the following fundamental identity that
we state as a lemma.

\begin{Lemma}
Assume that Hypotheses \ref{hp:introstateeq}, \ref{hp:introcost}, \ref
{hp:Hamiltonian1} and \ref{hp:exsolHJBvera} hold. Let $v\in C^{0,1}\left( %
\left[ 0,T\right] \times \mathbb{R}^{n}\right) $ be a strong solution of (%
\ref{eq:HJBristretta}). Then, for every $\left( t,x\right) \in \left[ 0,T%
\right] \times \mathbb{R}^{n}$ and $z\in \mathcal{Z}_{ad}\left( t\right) $
such that $J\left( t,x;z\right) <+\infty $ the following identity holds
\begin{equation*}
J\left( t,x;z\right) =v\left( t,x\right)
\end{equation*}
\begin{equation}
+\mathbb{E}\int_{t}^{T}\left[ -H^{0}\left( s,y(s),\partial _{x}v\left(
s,y\left( s\right) \right) \right) +H_{CV}^{0}\left( s,y(s),\partial
_{x}v\left( s,y\left( s\right) \right) ;z(s)\right) \right] ds
\label{eq:idfond1}
\end{equation}
where $y(s)\overset{def}{=}y\left( s;t,x,z\right) $ is the solution of (\ref
{eq:stateintro}) associated with the control $z$.
\end{Lemma}

\proof%
%
For the sake of completeness we first show how the proof goes in the case
when $v\in C^{1,2}\left( \left[ 0,T\right] \times \mathbb{R}^{n}\right) $
and is a strict solution of (\ref{eq:HJBristretta}). We use It\^{o}'s
formula applied to the function $v\left( t,x\right) $ by obtaining, for
every $\left( t,x\right) \in \left[ 0,T\right] \times \mathbb{R}^{n}$,
\begin{eqnarray}
v\left( T,y\left( T\right) \right) -v\left( t,y\left( t\right) \right)
&=&\int_{t}^{T}\left\langle \partial _{x}v\left( s,y\left( s\right) \right)
,B\left( s,y\left( s\right) \right) dW\left( s\right) \right\rangle  \notag
\\
&&+\int_{t}^{T}\left\langle F_{0}\left( s,y\left( s\right) \right) ,\partial
_{x}v\left( s,y\left( s\right) \right) \right\rangle ds
\label{eq:decversmooth} \\
&&+\int_{t}^{T}\left\langle F_{1}\left( s,y\left( s\right) ,z\left( s\right)
\right) ,\partial _{x}v\left( s,y\left( s\right) \right) \right\rangle ds
\notag \\
&&+\int_{t}^{T}\partial _{s}v\left( s,y\left( s\right) \right) ds  \notag \\
&&+\frac{1}{2}\int_{t}^{T}Tr\left[ B^{\ast }\left( s,y\left( s\right)
\right) \partial _{xx}v\left( s,y\left( s\right) \right) B\left( s,y\left(
s\right) \right) \right] ds,  \notag
\end{eqnarray}
(which is exactly the decomposition of Proposition 2.4 of \cite{FRFGpart1})
so, by taking expectation of both sides (this is finite since we have the
polynomial growth assumption) we get the so-called Dynkin formula
\begin{equation*}
\mathbb{E}\,v\left( T,y\left( T\right) \right) =\mathbb{E}\, v\left(
t,y\left( t\right) \right)
+\mathbb{E}\int_{t}^{T}\!\!\!\mathcal{L}_{0}v\left( s,y\left(
s\right) \right) ds+\mathbb{E}\int_{t}^{T}\!\!\!\left\langle
F_{1}\left( s,y\left( s\right) ,z\left( s\right) \right) ,\partial
_{x}v\left( s,y\left( s\right) \right) \right\rangle ds.
\end{equation*}
Now by (\ref{eq:HJBristretta}) we get, for every $s\geq 0$,
\begin{equation*}
\mathcal{L}_{0}v\left( s,y\left( s\right) \right) =-H^{0}\left( s,y\left(
s\right) ,\partial _{x}v\left( s,y\left( s\right) \right) \right) ,
\end{equation*}
which yields
\begin{equation*}
\mathbb{E}\,v\left( T,y\left( T\right) \right) =\mathbb{E}\,v\left(
t,y\left( t\right) \right)
\end{equation*}
\begin{equation*}
+\mathbb{E}\int_{t}^{T}\left[ -H^{0}\left( s,y\left( s\right) ,\partial
_{x}v\left( s,y\left( s\right) \right) \right) +\left\langle F_{1}\left(
s,y\left( s\right) ,z\left( s\right) \right) ,\partial _{x}v\left( s,y\left(
s\right) \right) \right\rangle \right] ds.
\end{equation*}
Now since we assumed that $J\left( t,x;z\right) <+\infty $ and since
Hypothesis \ref{hp:introcost} implies $J\left( t,x;z\right) >-\infty $ (this
also follows observing that the last line is not smaller than $-J\left(
t,x;z\right) $ and that the first line is finite by the polynomial growth of
$v$), we add $J\left( t,x;z\right) $ to both terms, use that $v\left(
T,y\left( T\right) \right) =\phi \left( y\left( T\right) \right) $ and that $%
y\left( t\right) =x$ to get
\begin{equation*}
J\left( t,x;z\right) =v\left( t,x\right) +\mathbb{E}\int_{t}^{T}-H^{0}\left(
s,y\left( s\right) ,\partial _{x}v\left( s,y\left( s\right) \right) \right)
ds
\end{equation*}
\begin{equation*}
+\mathbb{E}\int_{t}^{T}\left[ \left\langle F_{1}\left( s,y\left( s\right)
,z\left( s\right) \right) ,\partial _{x}v\left( s,y\left( s\right) \right)
\right\rangle +l\left( s,y\left( s\right) ,z\left( s\right) \right) \right]
ds
\end{equation*}
which is the claim (\ref{eq:idfond1}).

\bigskip

We now show how the proof goes when $v\in C^{0,1}\left( \left[ 0,T\right]
\times \mathbb{R}^{n}\right) $. The difference is due to the fact that the
term
\begin{equation*}
I:=\int_{t}^{T}\partial _{s}v\left( s,y\left( s\right) \right) ds+\frac{1}{2}%
\int_{t}^{T}Tr\left[ B^{\ast }\left( s,y\left( s\right) \right) \partial
_{xx}v\left( s,y\left( s\right) \right) B\left( s,y\left( s\right) \right) %
\right] ds
\end{equation*}
in the third and fourth line of equation (\ref{eq:decversmooth}) is not well
defined now. However by Hypothesis \ref{hp:exsolHJBvera} we know that $v$ is
a strong solution (in the sense of Definition \ref{df:solstrong}) of
\begin{equation*}
\mathcal{L}_{0}v\left( t,x\right) =h_{0}\left( t,x\right) ,\qquad v\left(
T,x\right) =\phi \left( x\right) ,
\end{equation*}
where we have set $h_{0}\left( t,x\right) =-H^{0}\left( t,x,\partial
_{x}v\left( t,x\right) \right) $. Then applying now Theorem \ref
{th:Itostrong} for the operator $\mathcal{L}_{0}$ (setting $b=F_{0}$, $%
\sigma =B$ and $b_{1}=F_{0}+F_{1}$) we get that
\begin{eqnarray*}
&&v\left( T,y\left( T\right) \right) =v\left( t,y\left( t\right) \right)
+\int_{t}^{T}\left\langle \partial _{x}v\left( s,y\left( s\right) \right)
,B\left( s,y\left( s\right) \right) dW\left( s\right) \right\rangle \\
&&+\int_{t}^{T}\left[ -H^{0}\left( s,y\left( s\right) ,\partial _{x}v\left(
s,y\left( s\right) \right) \right) +\left\langle F_{1}\left( s,y\left(
s\right) ,z\left( s\right) \right) ,\partial _{x}v\left( s,y\left( s\right)
\right) \right\rangle \right] ds.
\end{eqnarray*}
Now, taking expectation, adding and subtracting $J\left( t,x;z\right) $
(which is a.s. finite using the same argument for the smooth case), using
that $v\left( T,y\left( T\right) \right) =\phi \left( y\left( T\right)
\right) $ and that $y\left( t\right) =x$ we get
\begin{equation*}
\mathbb{E}\left[ \int_{t}^{T}l\left( s,y\left( s\right) ,z\left( s\right)
\right) ds+\phi \left( y\left( T\right) \right) \right] =v\left( t,x\right)
\end{equation*}
\begin{equation*}
+\mathbb{E}\int_{t}^{T}\left[ -H^{0}\left( s,y\left( s\right) ,\partial
_{x}v\left( s,y\left( s\right) \right) \right) +H_{CV}^{0}\left( s,y\left(
s\right) ,\partial _{x}v\left( s,y\left( s\right) \right) ;z\left( s\right)
\right) \right] ds
\end{equation*}
\smallskip and the claim again follows.\hfill
\hbox{\hskip 6pt\vrule width6pt height7pt
depth1pt  \hskip1pt}\bigskip

\textbf{Proof of Theorem \ref{th:VTintronostro}. }By the definition
of $H^{0} $ and $H_{CV}^{0}$, for every $\left( t,x\right) \in
\left[ 0,T\right] \times \mathbb{R}^{n}$, $z\in
\mathcal{Z}_{ad}\left( t\right) $ and for a.e. $\! \!\!s\in
\mathcal{T}_{t}$, the following inequality holds $\mathbb{P}$-almost
surely (setting $y\left( s\right) =y\left( s;t,x,z\right) $)
\begin{equation}
-H^{0}\left( s,y\left( s\right) ,\partial _{x}v\left( s,y\left( s\right)
\right) \right) +H_{CV}^{0}\left( s,y\left( s\right) ,\partial _{x}v\left(
s,y\left( s\right) \right) ;z\left( s\right) \right) \geq 0,
\label{eq:distriv}
\end{equation}
and then by the fundamental identity (\ref{eq:idfond1}) it follows $v\left(
t,x\right) \leq J\left( t,x;z\right) $ for every admissible $z$. This gives $%
v\leq V$ and so part (i) of the Theorem.

Now consider an admissible control $z$ such that, for a.e. \thinspace $s\in %
\left[ t,T\right] $, $\mathbb{P}-a.s.$%
\begin{equation}
H^{0}\left( s,y\left( s\right) ,\partial _{x}v\left( s,y\left( s\right)
\right) \right) =H_{CV}^{0}\left( s,y\left( s\right) ,\partial _{x}v\left(
s,y\left( s\right) \right) ;z\left( s\right) \right) .  \label{eq:scproof}
\end{equation}

By the fundamental identity (\ref{eq:idfond1}) we have
\begin{equation*}
v\left( t,x\right) =J\left( t,x;z\right) ,
\end{equation*}

which implies optimality and $v\left( t,x\right) =V\left( t,x\right) $%
.\hfill \hbox{\hskip 6pt\vrule width6pt height7pt
depth1pt  \hskip1pt}\bigskip

\begin{Remark}
\label{rm:novitasenzaderbis}

To get the methodological novelty in the above proof it is crucial to note
the following. Since we know that $u$ is the limit of classical solutions,
the first idea to approach the above proof (and that has been used e.g. in
the papers \cite{CPDE,EFJDE}), is probably to prove the fundamental identity
for the approximating solutions $u_{n}$ and then pass to the limit for $%
n\rightarrow \infty $. However if one tries to do this one needs to use the
uniform convergence of the space derivatives $\partial _{x}u_{n}$ to $%
\partial _{x}u$. This is not needed with our method, see also on this Remark
\ref{rm:suffperconvder}.\hfill
\hbox{\hskip 6pt\vrule
width6pt height7pt depth1pt  \hskip1pt}\bigskip
\end{Remark}

Let us go now to the proof of Verification Theorem \ref{th:VTintroinfhor}
for the infinite horizon case. Set $T=+\infty $, assume that $F_{0}$, $F_{1}$
and $B$ do not depend on time and that $l\left( t,x,z\right) =e^{\lambda
t}l_{1}\left( x,z\right) $ ($\lambda >0$), $\phi =0$; then we can set $t=0$
(since the dependence of the value function $V$ on $t$ becomes trivial in
this case) and the HJB equation for $V_{0}$ becomes an elliptic PDE.

We define the operator
\begin{equation*}
L_{0}:D\left( L_{0}\right) \subseteq C^{0}\left( \mathbb{R}^{n}\right)
\longrightarrow C^{0}\left( \mathbb{R}^{n}\right) ,\qquad D\left(
L_{0}\right) =C^{2}\left( \mathbb{R}^{n}\right) ,
\end{equation*}
\begin{equation*}
L_{0}u\left( x\right) =\left\langle F_{0}\left( x\right) ,\partial
_{x}u\left( x\right) \right\rangle +\frac{1}{2}\text{Tr }\left[ B^{\ast
}\left( x\right) \partial _{xx}u\left( x\right) B\left( x\right) \right]
\end{equation*}
and we rewrite the HJB equation (\ref{eq:HJBintroell}) as
\begin{equation}
\lambda u\left( x\right) +L_{0}u\left( x\right) +H^{1}\left( x,\partial
_{x}u\left( x\right) \right) =0,\quad x\in \mathbb{R}^{n}.  \label{eq:HJBell}
\end{equation}
where$\ H^{1}$ is defined in (\ref{eq:H1def}). We now define strong
solutions for such equations. Consider first, for $h\in C^{0}\left( \mathbb{R%
}^{n}\right) $ the inhomogenous elliptic problem
\begin{equation}
\lambda u\left( x\right) +L_{0}u\left( x\right) +h\left( x\right) =0,\quad
\quad \forall x\in \mathbb{R}^{n}.  \label{eq:elllinear}
\end{equation}

\begin{Definition}
We say that $u$ is a strict solution to the elliptic problem (\ref
{eq:elllinear}) if $u\in D\left( L_{0}\right) $ and (\ref{eq:elllinear})
holds.
\end{Definition}

\begin{Definition}
\label{df:solstrongell}We say that $u$ is a strong solution to the elliptic
problem (\ref{eq:elllinear}) if there exists a sequence $\left( u_{n}\right)
\subseteq D\left( L_{0}\right) $ and a sequence $\left( h_{n}\right)
\subseteq C^{0}\left( \mathbb{R}^{n}\right) $, such that

\begin{enumerate}
\item  For every $n\in \mathbb{N}$ $u_{n}$ is a strict solution of the
problem
\begin{equation*}
\lambda u_{n}\left( x\right) -L_{0}u_{n}\left( x\right) =h_{n}\left(
x\right) ,\qquad \quad \quad \forall x\in \mathbb{R}^{n}.
\end{equation*}

\item  The following limits hold
\begin{equation*}
\begin{array}{l}
u_{n}\longrightarrow u\text{ in }C^{0}\left( \mathbb{R}^{n}\right) , \\
h_{n}\longrightarrow h\text{ in }C^{0}\left( \mathbb{R}^{n}\right) .
\end{array}
\end{equation*}
\end{enumerate}
\end{Definition}

For functions $u\in C^{1}\left( \mathbb{R}^{n}\right) $ that are strong
solutions of the Cauchy problem (\ref{eq:elllinear}) the result of Theorem
\ref{th:Itostrong} still holds: it is indeed a simpler case. Now we go to
strong solutions of the HJB equation (\ref{eq:HJBell}).

\begin{Definition}
\label{df:solHJBstrongell}A function $v\in C^{1}\left( \mathbb{R}^{n}\right)
$ is a strong solution of the HJB equation (\ref{eq:HJBell}) if, setting $%
h_{0}\left( x\right) =-H^{1}\left( x,\partial _{x}v\left( x\right) \right) $%
, $v$ is a strong solution of the linear problem
\begin{equation*}
\lambda v\left( x\right) +L_{0}v\left( x\right) =h_{0}\left( x\right) ,
\end{equation*}
in the sense of Definition \ref{df:solstrongell}.
\end{Definition}

Assume the following.

\begin{Hypothesis}
\label{hp:exsolHJBverabis}There exists a function $v\in C^{1}\left( \mathbb{R%
}^{n}\right) $ which is a strong solution of the HJB equation (\ref
{eq:HJBell}) in the sense of Definition \ref{df:solHJBstrongell}. Moreover:

\begin{itemize}
\item[(i)]  either we can choose the approximating sequence $\left(
v_{n}\right) $ of Definition \ref{df:solHJBstrongell} so that for every $%
s\geq 0$ and for every admissible control $z\in \mathcal{Z}_{ad}\left(
0\right) $%
\begin{equation*}
\lim_{n\rightarrow +\infty }\int_{0}^{s}\left\langle \partial
_{x}v_{n}\left( y\left( r\right) \right) -\partial _{x}v\left( y\left(
r\right) \right) ,F_{1}\left( y\left( r\right) ,z\left( r\right) \right)
\right\rangle dr=0,\qquad u.c.p.,
\end{equation*}

\item[(ii)]  or the function
\begin{equation*}
\left( x,z\right) \rightarrow B^{-1}\left( x\right) F_{1}\left( x,z\right) ,
\end{equation*}
where $B^{-1}$ stands for the pseudo-inverse of $B$, is well defined and
bounded on $\mathbb{R}^{n}\times U$.
\end{itemize}
\end{Hypothesis}

The precise statement of the verification theorem is the following.

\begin{Theorem}
\label{th:VTinfhor}Assume that Hypotheses \ref{hp:introstateeq}, and \ref
{hp:exsolHJBverabis} hold. Assume also that $l_{1}$ is continuous and
bounded and that Hypothesis \ref{hp:Hamiltonian1} hold true with $H^{1}$ in
place of $H^{0}$. If $v$\ is a bounded strong solution of the HJB equation (%
\ref{eq:HJBell}) and $v\in C^{1}\left( \mathbb{R}^{n}\right) $\ then $v\leq
V_{0}$ on $\mathbb{R}^{n}$. Morerover if $z$ is an admissible control at $%
\left( 0,x\right) $ that satisfies (setting $y\left( s\right) =y\left(
s;0,x,z\right) $)
\begin{equation*}
H^{1}\left( y\left( s\right) ,\partial _{x}v\left( y\left( s\right) \right)
\right) =H_{CV}^{1}\left( y\left( s\right) ,\partial _{x}v\left( y\left(
s\right) \right) ;z\left( s\right) \right)
\end{equation*}
for $a.e.$ $s\in \left[ 0,+\infty \right) $, $\mathbb{P-}$almost surely,
then $z$ is optimal and $v\left( x\right) =V_{0}\left( x\right) .$
\end{Theorem}

\textbf{Proof. }The proof goes along the same lines as the finite horizon
case. First we observe that the boundedness of the datum $l_{1}$ implies
that also $V_{0}$ is bounded and that for every $x\in \mathbb{R}^{n}$, $z\in
\mathcal{Z}_{ad}\left( 0\right) $ the functional $J\left( 0,x;z\right) $ is
finite. Then we prove the following fundamental identity
\begin{eqnarray*}
J\left( 0,x;z\right) =v\left( x\right) +\mathbb{E}\int_{0}^{+\infty
}\!\!\!\! \! e^{-\lambda s}\left[ -H^{1}\left( y(s),\partial _{x}v\left(
y\left( s\right) \right) \right) +H_{CV}^{1}\left( y(s),\partial _{x}v\left(
y\left( s\right) \right) ;z(s)\right) \right] ds.
\end{eqnarray*}
To do it we apply Theorem \ref{th:Itostrong} for the operator $\mathcal{L}%
_{0}=\partial _{t}+L_{0}$ (so $b=F_{0}$, $\sigma =B$ and $b_{1}=F_{0}+F_{1}$%
) but taking $t=0$ and replacing $v\left( s,x\right) $ by
$e^{-\lambda s}v\left( x\right) $ for $s\ge 0$. Observe first that,
if $v\in C^{2}\left( \mathbb{R}^{n}\right) $ then the function
$w\left( t,x\right) =e^{-\lambda t}v\left( x\right) $ solves, for
$t\ge 0$, $x \in \mathbb{R}^{n}$  the equation
\begin{equation}
w_{t}\left( t,x\right) +L_{0}w\left( t,x\right) =-e^{-\lambda t}H^{1}\left(
x,\partial _{x}v\left( x\right) \right) .  \label{eq:paracc}
\end{equation}
Moreover if $v\in C^{1}\left( \mathbb{R}^{n}\right) $ is a strong solution
of (\ref{eq:HJBell}), also $w\left( t,x\right) $ is a strong solution of (%
\ref{eq:paracc}) and satisfies the assumptions of Theorem \ref{th:Itostrong}%
. So we get for every $T_{1}>0$%
\begin{eqnarray*}
&&e^{-\lambda T_{1}}v\left( y\left( T_{1}\right) \right) =v\left( y\left(
0\right) \right) +\int_{0}^{T_{1}}e^{-\lambda s}\left\langle \partial
_{x}v\left( y\left( s\right) \right) ,B\left( y\left( s\right) \right)
dW\left( s\right) \right\rangle \\
&&+\int_{0}^{T_{1}}e^{-\lambda s}\left[ -H^{1}\left( s,y\left( s\right)
,\partial _{x}v\left( s,y\left( s\right) \right) \right) +\left\langle
F_{1}\left( s,y\left( s\right) ,z\left( s\right) \right) ,\partial
_{x}v\left( s,y\left( s\right) \right) \right\rangle \right] ds.
\end{eqnarray*}
Now, taking expectation, adding and subtracting $\mathbb{E}%
\int_{0}^{T_{1}}e^{-\lambda s}l_{1}\left( y\left( s\right) ,z\left( s\right)
\right) ds$ (which is a.s. finite by the boundedness of $l_{1}$) and using
that $y\left( 0\right) =x$ we get
\begin{eqnarray*}
&&\mathbb{E}\left[ \int_{0}^{T_{1}}e^{-\lambda s}l_{1}\left( y\left(
s\right) ,z\left( s\right) \right) ds+e^{-\lambda T_{1}}v\left( y\left(
T_{1}\right) \right) \right] \\
&=&v\left( x\right) +\mathbb{E}\int_{0}^{T_{1}}e^{-\lambda s}\left[
-H^{1}\left( y\left( s\right) ,\partial _{x}v\left( y\left( s\right) \right)
\right) +H_{CV}^{1}\left( y\left( s\right) ,\partial _{x}v\left( y\left(
s\right) \right) ;z\left( s\right) \right) \right] ds
\end{eqnarray*}
\smallskip and the claim again follows taking the limit for $%
T_{1}\rightarrow +\infty $ and using the boundedness of $v$. The rest of the
proof is exactly the same as in the finite horizon case.\hfill
\hbox{\hskip 6pt\vrule width6pt height7pt
depth1pt  \hskip1pt}\bigskip

\begin{Remark}
\label{rm:verinfhornobdd}It must be noted that the assumption on the
boundedness of $l_{1}$ and $v$ is made for simplicity of exposition. What is
really needed in this infinite horizon case is that $J$ is well defined for
each $z\in \mathcal{Z}_{ad}\left( 0\right) $ and that $\lim_{T_{1}%
\rightarrow +\infty }$ $e^{-\lambda T_{1}}v\left( y\left( T_{1}\right)
\right) =0$ for every $z\in \mathcal{Z}_{ad}\left( 0\right) $. This allows
to pass to the limit in the latter formula. In some cases this can be checked
directly, in other cases much weaker conditions can be imposed (e.g.
sublinearity of $l_{1}$ and estimates on the solution $y$ as $%
T_{1}\rightarrow +\infty $).\hfill
\hbox{\hskip 6pt\vrule width6pt height7pt
depth1pt  \hskip1pt}\bigskip
\end{Remark}

\section{Necessary Conditions and Optimal Feedback Controls\label{FEEDBACK}}

Here we want simply to show that under additional assumptions (mainly about
existence and/or uniqueness of the maximum of the Hamiltonian and of the
solution of the closed loop equation) we can get necessary conditions and
optimal feedback controls. This is a consequence of the verification
Theorems \ref{th:VTintronostro} and \ref{th:VTinfhor} which has some
importance for applications. For the sake of brevity we give the results
only for the finite horizon case observing that completely analogous results
hold true for the infinite horizon case. We start by making a remark about
the so-called ``weak formulation'' of a stochastic control problem.

\begin{Remark}
\label{rm:weakformulation}\textrm{The setting introduced in Section \ref
{INTROVER} corresponds to the so-called ``strong formulation'' of a
stochastic control problem. For certain purposes (namely the existence of
optimal feedbacks) it is convenient to consider the ``weak formulation'',
letting the stochastic basis vary. In such formulation one considers as the
set }$\mathcal{\overline{Z}}_{ad}\left( t\right) $\textrm{\ of admissible
controls as the set of 5-tuples }$\left( \Omega ,\mathcal{F},\mathbb{P}%
,W,z\right) $\textrm{\ such that}

\begin{itemize}
\item  $\left( \Omega ,\mathcal{F},\mathbb{P}\right) $\textrm{\ is a
complete filtered probability space with the filtration }$\mathcal{F}$%
\textrm{\ satisfying the usual conditions,}

\item  $W$\textrm{\ is an }$E$\textrm{\ valued }$m$\textrm{-dimensional
Brownian motion on }$\left[ t,T\right] $\textrm{,}

\item  \textrm{the process }$z$\textrm{\ is measurable, }$\left( \mathrm{%
\mathcal{F}_{s}}\right) $\textrm{-adapted and a.s. locally integrable.}
\end{itemize}

\textrm{We will use the notation }$\left( \Omega ,\mathcal{F},\mathbb{P}%
,W,z\right) \in \mathcal{\overline{Z}}_{ad}\left( t\right) $ \textrm{and we
will call this control strategies weakly admissible (weakly optimal when
they are optimal). When no ambiguity arises we will leave aside the
probability space and the white noise (regarding it as fixed) and simply
consider }$z\in \mathcal{Z}_{ad}\left( t\right) $\textrm{. We will specify
when the weak concept of admissible control is needed. The same will be done
in Sections \ref{APPL1EXIT}\ and \ref{COMPARISON}. One can see e.g.
\cite[p.64]{Zhoubook}, \cite[pp.141, 160]{FS} for comments on these
formulations.\hfill \hfill
\hbox{\hskip 6pt\vrule width6pt height7pt
depth1pt  \hskip1pt}\bigskip }
\end{Remark}

We start recalling a case when the sufficient condition of the verification
Theorem \ref{th:VTintronostro} becomes also necessary.

\begin{Proposition}
\label{pr:feedback2}Assume that the hypotheses of Theorem \ref
{th:VTintronostro} hold true. We also assume that the value function $V$ is
a strong solution of the HJB equation. Then an admissible control $z\in
\mathcal{Z}_{ad}\left( t\right) $ is optimal at $\left( t,x\right) \in \left[
0,T\right] \times \mathbb{R}^{n}$ if and only if (\ref{eq:condsuff}) holds
with $V$ in place of $v$.
\end{Proposition}

\textbf{Proof.} The fundamental identity in this case gives
\begin{equation*}
J\left( t,x;z\right) =V\left( t,x\right)
\end{equation*}
\begin{equation*}
+\mathbb{E}\int_{t}^{T}\left[ H^{0}\left( s,y\left( s\right) ,\partial
_{x}V\left( s,y\left( s\right) \right) \right) -H_{CV}^{0}\left( s,y\left(
s\right) ,\partial _{x}V\left( s,y\left( s\right) \right) ;z\left( s\right)
\right) \right] ds
\end{equation*}
and this immediately gives the claim.\textrm{\hfill
\hbox{\hskip 6pt\vrule width6pt height7pt
depth1pt  \hskip1pt}\bigskip }

\begin{Remark}
\textrm{One case where the above Proposition \ref{pr:feedback2} can be
applied is when:}

\begin{itemize}
\item  \textrm{it is known that }$V$\textrm{\ is the unique viscosity
solution of the HJB equation;}

\item  \textrm{it is known that strong solutions are also viscosity
solutions.}
\end{itemize}

\textrm{This happens e.g. in the example of Section \ref{APPL2DEG}.\hfill
\hbox{\hskip 6pt\vrule width6pt height7pt
depth1pt  \hskip1pt}\bigskip }
\end{Remark}

We now pass to the existence of feedbacks. First we recall the definition of
(weakly) admissible feedback map.

\begin{Definition}
\label{df:feedback} A measurable map $G:\left[ 0,T\right] \times \mathbb{R}%
^{n}\longmapsto U$ is a (weakly) admissible feedback map if for every $%
\left( t,x\right) \in \left[ 0,T\right] \times \mathbb{R}^{n}$ the closed
loop equation ($s\in \left[ t,T\right] $)
\begin{equation}
y(s)=x+\int_{t}^{s}F_{1}\left( r,y\left( r\right) ,G\left( r,y\left(
r\right) \right) \right) dr+\int_{t}^{s}F_{0}\left( r,y\left( r\right)
\right) dr+\int_{t}^{s}B\left( r,y\left( r\right) \right) dW(r),
\label{eq:cleadm}
\end{equation}
admits a (weak) solution $y = y_{G,t,x}  $ for every $\left(
t,x\right) \in \left[ 0,T\right] \times \mathbb{R}^{n}$ and the control
strategy $z_{G,t,x}\left( s\right) =G\left( s,y_{G,t,x}\left( s\right)
\right) $ is (weakly) admissible.
\end{Definition}

If we have a (weakly) admissible feedback map $G$ then the control strategy $%
z_{G,t,x}\left( s\right) =G\left( s,y_{G,t,x}\left( s\right) \right) $ is by
definition (weakly) admissible and so, given $G$, we have, for every $\left(
t,x\right) \in \left[ 0,T\right] \times \mathbb{R}^{n}$ an admissible (weak)
couple $\left( z_{G,t,x},  y_{G,t,x}
\right) $. Our goal is to find an admissible feedback map $G$ such that, for
every $\left( t,x\right) \in \left[ 0,T\right] \times \mathbb{R}^{n}$, the
strategy $z_{G,t,x}$ is (weakly) optimal. Such $G$ will be called a (weakly)
optimal feedback map.

\begin{Proposition}
\label{pr:feedbackvero} Assume that the Hypotheses of Theorem \ref
{th:VTintronostro} hold true. Assume also that:

\begin{enumerate}
\item[(i)]  the maximum of the Hamiltonian (\ref{eq:H0def}) exists and it is
possible to define a measurable map
\begin{equation*}
G_{0}:\left[ 0,T\right] \times \mathbb{R}^{n}\times \mathbb{R}%
^{n}\longmapsto U
\end{equation*}
such that $G_{0}\left( t,x,p\right) \in \arg \min_{z}H_{CV}\left(
t,x,p;z\right) $ for every $\left( t,x,p\right) \in \left[ 0,T\right] \times
\mathbb{R}^{n}\times \mathbb{R}^{n}$;

\item[(ii)]  the map $G\left( t,x\right) =G_{0}\left( t,x,\partial
_{x}v\left( t,x\right) \right) $ is a (weakly) admissible feedback map.
\end{enumerate}

Then for every $\left( t,x\right) \in \left[ 0,T\right] \times \mathbb{R}%
^{n} $, $z_{G,t,x}$ is a (weakly) optimal control, $y_{G,t,x}$ is the
corresponding (weakly) optimal state and $v=V$.

Finally if $\arg \min_{z\in U}H_{CV}\left( t,x,v_{x}(t,x);z\right) $ is
always a singleton and the closed loop equation (\ref{eq:cleadm}) admits a
unique (weak) solution, then the optimal control is (weakly) unique.
\end{Proposition}

\textbf{Proof. }This is an obvious consequence of the verification Theorem
\ref{th:VTintronostro}.\hfill
\hbox{\hskip 6pt\vrule width6pt height7pt
depth1pt  \hskip1pt}\bigskip

\begin{Remark}
We observe that if we have existence of the (weak) solution of the closed
loop equation (\ref{eq:cleadm}) for every strong solution of the HJB
equation (\ref{eq:HJBristretta}), then we automatically have that the strong
solution is unique. \hfill
\hbox{\hskip 6pt\vrule width6pt height7pt
depth1pt  \hskip1pt}\bigskip
\end{Remark}

\begin{Remark}
We finally observe that Propositions \ref{pr:feedback2} and \ref
{pr:feedbackvero} holds also in the infinite horizon case with obvious
changes. We will use the infinite horizon version of Proposition \ref
{pr:feedbackvero} in the study of the one dimensional example in Section \ref
{APPL2DEG}.\hfill
\hbox{\hskip
6pt\vrule width6pt height7pt depth1pt  \hskip1pt}\bigskip
\end{Remark}

\section{Application 1: a class of exit time control problems with non
degenerate diffusion\label{APPL1EXIT}}

Here we show that our technique for proving verification theorems works for
a family of stochastic optimal control problems with exit time and
nondegenerate diffusion where the solutions of the associated HJB equation
are not known to be $C^{1,2}$.

\subsection{The problem}

Let first $\mathbb{R}^{n}$ be the state space, $\mathbb{R}^{n}$ be the space
of noises and $U$ (a given subset of a Polish space) be the control space.
Let $T$ be a fixed finite horizon, $\left( \Omega ,\mathcal{F},\mathbb{P}%
\right) $ be a given stochastic basis (where $\mathcal{F}$ stands for a
given filtration $\left( \mathcal{F}_{s}\right) _{s\in \left[ 0,T\right] }$
satisfying the usual conditions), $W$ be a cylindrical Brownian motion with
values in $\mathbb{R}^{n}$ adapted to $\left( \mathcal{F}_{s}\right) _{s\in %
\left[ 0,T\right] }$. Consider a stochastic controlled system in $\mathbb{R}%
^{n}$ with fixed finite horizon $T$ and initial time $t\in \left[ 0,T\right)
$ governed by the state equation
\begin{equation}
\left\{
\begin{array}{ll}
dy\left( s\right) =\left[ F_{0}\left( s,y\left( s\right) \right)
+F_{1}\left( s,y\left( s\right) ,z\left( s\right) \right) \right] ds+B\left(
y\left( s\right) \right) dW(s),\qquad s\in \left[ t,T\right]  &  \\
y\left( t\right) =x. &
\end{array}
\right.   \label{eq:SE1}
\end{equation}
We consider an \emph{open bounded domain} $\mathcal{O}\subseteq \mathbb{R}%
^{n}$ with uniformly $C^{2}$ boundary (see e.g. \cite[pp.2-3]{Lunardibook}
for the definition). The initial datum $x$ belongs to $\mathcal{O}$ and we
call $\tau _{\mathcal{O}}$ the first exit time of the process $y$ from this
open set, i.e.
\begin{equation*}
\tau _{\mathcal{O}}\left( \omega \right) =\inf \left\{ s>t:y\left( s;\omega
\right) \in \mathcal{O}^{c}\right\} .
\end{equation*}
Moreover $F_{0}:\left[ 0,T\right] \times \bar{\mathcal{O}}\rightarrow
\mathbb{R}^{n}$, $F_{1}:\left[ 0,T\right] \times \bar{\mathcal{O}}\times
U\rightarrow \mathbb{R}^{n}$, $B:\left[ 0,T\right] \times \bar{\mathcal{O}}%
\rightarrow \mathcal{L}\left( \mathbb{R}^{n};\mathbb{R}^{n}\right) $, ($\bar{%
\mathcal{O}}$ stands for the closure of the set $\mathcal{O}$) satisfy

\begin{Hypothesis}
\label{hp:appl1driftbounded}$F_{0}$ and $F_{1}$ are continuous and bounded.
\end{Hypothesis}

\begin{Hypothesis}
\label{hp:appl1nondeg}$B$ is time independent, uniformly continuous and
bounded and nondegenerate i.e. there exists $\lambda _{0}>0$ such that
\begin{equation*}
\lambda _{0}^{-1}\left| \left| \xi \right| \right| \geq \sum_{i,j=1}^{n}
\left[ B\left( x\right) B^{T}\left( x\right) \right] _{i,j}\xi _{i}\xi
_{j}\geq \lambda _{0}\left| \left| \xi \right| \right| \qquad \forall \left(
t,x\right) \in \left[ 0,T\right] \times \bar{\mathcal{O}},\qquad \forall \xi
\in \mathbb{R}^{n}.
\end{equation*}
\end{Hypothesis}

\begin{Remark}
First note that the above Hypothesis \ref{hp:appl1nondeg} implies that $B$
is invertible with bounded inverse, so we fall in the Hypothesis \ref
{hp:exsolHJBvera}-(ii). Moreover it is possible to extend our results to the
case when $B$ is time dependent, provided suitable regularity conditions are
satisfied. Such conditions are needed to apply the results of semigroup
theory in the subsection below and substantially require that the domain of
the operator
\begin{equation*}
L\left( t\right) u=\frac{1}{2}\mathrm{Tr}\left[ B^{\ast }\left( t,x\right)
B\left( t,x\right) \partial _{xx}u\right]
\end{equation*}
be constant and the coefficient $B$ be H\"{o}lder continuous in time,
see on this \cite{Buttu1,Buttu2} and \cite[Section 3.1.2]{Lunardibook}. We
leave aside these assumptions to keep a more simplified setting where,
anyway, $C^{1,2}$ regularity does not hold in general.\textrm{\hfill
\hbox{\hskip 6pt\vrule width6pt height7pt depth1pt
\hskip1pt}\bigskip }
\end{Remark}

About the control strategy $z:[t,T]\times \Omega \rightarrow U$ we assume
that it belongs to a given set $\mathcal{Z}_{ad}\left( t\right) $ of
stochastic processes defined on $\left[ t,T\right] \times \Omega $ with
values in a fixed Polish space $U$. More precisely we will assume the
following.

\begin{Hypothesis}
\label{hp:appl1control}The control space $U$ is a Polish space, while $%
\mathcal{Z}_{ad}\left( t\right) $ is the space of all measurable processes $%
z:[t,T]\times \Omega \rightarrow U$ adapted to the filtration $\mathcal{F}$.
\end{Hypothesis}

\begin{Remark}
\textrm{The above setting corresponds to the so-called ``strong
formulation'' of a stochastic control problem. In this problem, since the
state equation admits only weak solutions in general (see next proposition)
we need to consider the ``weak formulation'', letting the stochastic basis
to vary, see Remark \ref{rm:weakformulation}. To avoid heavy notation we
will keep the same symbols as in the strong formulation, as we did in
Section \ref{FEEDBACK}. \hfill
\hbox{\hskip 6pt\vrule width6pt height7pt depth1pt
\hskip1pt}\bigskip }
\end{Remark}


The following Proposition can be proved in the same way as in \cite{SV},
even if the drift is random. In fact Theorem 6.4.3 therein makes use of
Girsanov transformation to eliminate the drift; the same can be done here.
Then we can apply Corollary 6.4.4 and Theorem 7.2.1 of \cite{SV}.

\begin{Proposition}
Assume that Hypotheses \ref{hp:appl1driftbounded}, \ref{hp:appl1nondeg}
hold. Then, for all $z\in \mathcal{Z}_{ad}\left( t\right) $, equation (\ref
{eq:SE1}) has a weak solution
\begin{equation*}
y(\cdot ;t,x,z)\in C_{\mathcal{F}}^{0}\left( t,T;X\right) .
\end{equation*}
This solution is unique in the sense of probability law.
\end{Proposition}

\begin{Remark}
\textrm{In the case of dimension 1, no continuity on the diffusion
coefficients is required, see \cite{SV} Exercise 7.3.3 at page 192.\hfill
\hbox{\hskip 6pt\vrule width6pt height7pt
depth1pt  \hskip1pt}\bigskip }
\end{Remark}

We now consider the following stochastic optimal control problem with exit
time. Minimize the cost functional
\begin{eqnarray}
J\left( t,x;z\right)  &=&\mathbb{E}\left[ \int_{t}^{\tau _{\mathcal{O}%
}\wedge T}l\left( s,y\left( s;t,x,z\right) ,z\left( s\right) \right)
ds\right.   \label{eq:CFappl1} \\[0.1in]
&&\left. +I_{\left\{ \tau _{\mathcal{O}}<T\right\} }\psi \left( \tau _{%
\mathcal{O}},y\left( \tau _{\mathcal{O}};t,x,z\right) \right) +I_{\left\{
\tau _{\mathcal{O}}\geq T\right\} }\phi \left( y\left( T;t,x,z\right)
\right) \right] ,  \notag
\end{eqnarray}
over all controls $z\in \mathcal{Z}_{ad}\left( t\right) $. Here $y\left(
\cdot ;t,x,z\right) $ is the solution of the equation (\ref{eq:SE1}) and we
assume that $l,\psi ,\phi $ satisfy

\begin{Hypothesis}
\label{hp:appl1functional}$l\in C^{0}\left( \left[ 0,T\right] \times \bar{%
\mathcal{O}}\right) $, $\phi \in C^{0}\left( \bar{\mathcal{O}}\right) $ $%
\psi \in C^{\frac{1+\beta }{2},1+\beta }\left( \left[ 0,T\right] \times
\mathbb{\partial }\mathcal{O}\right) $ (for some $\beta >0$), $\psi \left(
T,x\right) =\phi \left( x\right) $ on $\partial \mathcal{O}$.
\end{Hypothesis}

\begin{Remark}
The above hypothesis is needed to apply the results of semigroup theory in
the subsection below. In particular, under Hypothesis \ref
{hp:appl1functional} the operator
\begin{equation*}
Lu=\frac{1}{2}\mathrm{Tr}\left[ B^{\ast }\left( x\right) B\left( x\right)
\partial _{xx}u\right]
\end{equation*}
(with 0 Dirichlet boundary conditions) generates an analytic semigroup and
the associated boundary value problem in $\mathcal{O}$ is well posed (see
\cite[Ch.5]{Lunardibook}).\hfill \textrm{%
\hbox{\hskip 6pt\vrule width6pt height7pt
depth1pt  \hskip1pt}\bigskip }
\end{Remark}

The value function of this problem is defined as
\begin{equation}
V(t,x)=\inf \left\{ J\left( t,x;z\right) \,:\;z\in \mathcal{Z}_{ad}\left(
t\right) \right\}   \label{eq:VFappl1}
\end{equation}
and a control $z^{\ast }\in \mathcal{Z}_{ad}\left( t\right) $ and such that $%
V(t,x)=J(t,x;z^{\ast })$ is said to be \emph{optimal} with respect to the
initial time and state $\left( t,x\right) $. The corresponding HJB equation
is
\begin{equation}  \label{eq:HJBappl1}
\left\{
\begin{array}{l}
v\left( t,x\right) +\frac{1}{2}\mathrm{Tr}\left[ B\left( x\right) \partial
_{xx}v(t,x)B^{\ast }\left( x\right) \right] =\left\langle F_{0}\left(
t,x\right) ,\partial _{x}v\left( t,x\right) \right\rangle +H^{0}\left(
t,x,\partial _{x}v\left( t,x\right) \right) , \\
\qquad \qquad \qquad \qquad \qquad \qquad \qquad \qquad \qquad \qquad \qquad
t\in \left[ 0,T\right] \text{,\qquad }x\in \mathcal{O}, \\
v\left( T,x\right) =\phi \left( x\right) ,\qquad x\in \mathcal{O}, \\
v\left( t,x\right) =\psi \left( t,x\right) ,\quad \quad t\in \left[ 0,T%
\right] ,x\in \partial \mathcal{O},
\end{array}
\right.
\end{equation}
where
\begin{equation}
H^{0}\left( t,x,p\right) =\inf_{z\in U}H_{CV}^{0}(t,x,p;z),
\label{eq:appl1Hamiltonian}
\end{equation}
with
\begin{equation*}
H_{CV}^{0}(t,x,p;z)=\left\langle F_{1}\left( t,x,z\right) ,p\right\rangle
+l\left( t,x,z\right) ,
\end{equation*}
being the \textit{current value Hamiltonian}.

\begin{Proposition}
\label{pr:regHamappl1}Under Hypotheses \ref{hp:appl1driftbounded}, \ref
{hp:appl1nondeg}, \ref{hp:appl1control}, \ref{hp:appl1functional}, the
Hamiltonian $H^{0}\left( t,x,p\right) $ defined in (\ref{eq:appl1Hamiltonian}%
) is continuous in $\left[ 0,T\right] \times \bar{\mathcal{O}}\times \mathbb{%
R}^{n}$. Moreover there exists a constant $C>0$ such that
\begin{eqnarray}
\left| H^{0}\left( t,x,p\right) -H^{0}\left( t,x,q\right) \right|  &\leq
&C\left| p-q\right|   \notag \\
\left| H^{0}\left( t,x,p\right) \right|  &\leq &C\left( 1+\left| p\right|
\right)   \label{eq:lipHamappl1}
\end{eqnarray}
\end{Proposition}

\textbf{Proof.} It is enough to apply the definition of the Hamiltonian and
use the continuity and the boundedness of $F_1$ and $l$. \hfill \textrm{%
\hbox{\hskip 6pt\vrule width6pt height7pt depth1pt
\hskip1pt}\bigskip }

\subsection{Strong solutions of the HJB equation\label{STRONGHJB}}

The HJB equation (\ref{eq:HJBappl1}) above is a semilinear parabolic
equation with continuous coefficients. Since the second order term is
nondegenerate (Hypothesis \ref{hp:appl1nondeg}) one expects interior
regularity results for the solution even if the boundary data are merely
continuous, as it is under our assumptions. In particular the following
result, taken from \cite{CKSS}, Theorems 9.1 and 9.2, apply.

\begin{Theorem}
\label{th:Swiech}Under Hypotheses \ref{hp:appl1driftbounded}, \ref
{hp:appl1nondeg}, \ref{hp:appl1control}, \ref{hp:appl1functional}, there
exists a viscosity solution $u$ of (\ref{eq:HJBappl1}) and $u$ belongs to
the space $C_{loc}^{\frac{1+\alpha }{2},1+\alpha }\left( \left[ 0,T\right)
\times \mathcal{O}\right) \cap C^{0}\left( \left[ 0,T\right] \times \bar{%
\mathcal{O}}\right) $ for some $\alpha \in \left( 0,1\right) $.
\end{Theorem}

We note also that it is not known if the solution is classical in the
interior. This is known, using theorems of \cite{caffarelli}, if the data
are supposed to be H\"{o}lder continuous in $\left( t,x\right) $, which we
do not assume here. Moreover even in the linear case, if the coefficients
are only continuous, solutions are only $W^{2,p}$, so they are also $%
C^{1,\alpha }$. An example in this direction is in \cite[Ch.4]{GT}.

This means that we are exactly in the case where it makes sense to apply our
technique to prove verification theorems: no classical solution but
solutions (in a generalized sense) with at least $C^{0,1}$ regularity. Once
this is clear, we need to check if our solutions (that exists at least in
the viscosity sense thanks to the above Theorem \ref{th:Swiech}) are also
strong solutions in the sense of Definition \ref{df:solHJBstrong} (suitably
modified to take care of the boundary datum $\psi $). Such a result is not
available in the literature in this form but it can be easily deduced using
the results on analytic semigroups contained e.g. in \cite{Lunardibook}. We
explain here below how this can be done in the case when $\psi =0$ (the
general case can be treated with the ideas explained in \cite[Section 5.1.2]
{Lunardibook}). In this case the operator
\begin{eqnarray*}
L &:&D\left( L\right) \subseteq C^{0}\left( \bar{\mathcal{O}}\right)
\rightarrow C^{0}\left( \bar{\mathcal{O}}\right) , \\
D\left( L\right)  &=&\left\{ \eta \in \cap _{p\geq 1}W_{loc}^{2,p}\left(
\mathcal{O}\right) :\quad \eta ,L\eta \in C^{0}\left( \bar{\mathcal{O}}%
\right) ,\quad \eta |_{\partial \mathcal{O}}=0\right\} , \\
\left( L\eta \right) \left( x\right)  &=&\frac{1}{2}\mathrm{Tr}\left[
B\left( x\right) \partial _{xx}\eta (x)B^{\ast }\left( x\right) \right] ,
\end{eqnarray*}
generates an analytic semigroup $\left\{ e^{tL},t\geq 0\right\} $,
\cite[p.97, Corollary 3.1.21]{Lunardibook}. Moreover, given any initial
datum $\phi \in C^{0}\left( \bar{\mathcal{O}}\right) $ and any function $%
H^{0}$ continuous in $\left[ 0,T\right] \times \bar{\mathcal{O}}\times
\mathbb{R}^{n}$ and satisfying (\ref{eq:lipHamappl1}) we can apply a
modification of the argument used to prove Proposition 7.3.4 in \cite[p.281]
{Lunardibook} to get existence and uniqueness of a solution $u$ of the
integral equation
\begin{eqnarray}
u\left( t,x\right)  &=&\left( e^{\left( T-t\right) L}\phi \right) \left(
x\right)   \label{eq:HJBmildappl1} \\
&&+\int_{t}^{T}\left( e^{\left( T-s\right) L}\left[ \left\langle F_{0}\left(
s,\cdot \right) ,\partial _{x}u\left( s,\cdot \right) \right\rangle
+H^{0}\left( s,\cdot ,\partial _{x}u\left( s,\cdot \right) \right) \right]
\right) \left( x\right) ds,  \notag
\end{eqnarray}
which can be considered as an integral form of the PDE (\ref{eq:HJBappl1})
when $\psi =0$ (in the general case one needs to lift the function $\psi $
into the equation obtaining an extra term in (\ref{eq:HJBmildappl1}), see on
this \cite[Remark 5.1.14, p.195]{Lunardibook}) written using the variation
of constants formula. More precisely we have the following definitions and
results.

\begin{Definition}
We say that a function $w$ belongs to the space $\Sigma ^{1,\alpha }\left( %
\left[ 0,T\right] \times \bar{\mathcal{O}}\right) $ if $w\in C^{0}\left( %
\left[ 0,T\right] \times \bar{\mathcal{O}}\right) $, $w$ is Fr\'{e}chet
differentiable in $x\in \mathcal{O}$, $\partial _{x}w\in C^{0}\left( \left[
0,T-\varepsilon \right] \times \bar{\mathcal{O}}\right) $ for every $%
\varepsilon \in \left( 0,T\right) $, and
\begin{equation*}
\sup_{\left( t,x\right) \in \left[ 0,T\right) \times \bar{\mathcal{O}}%
}\left| \left( T-t\right) ^{\alpha }\partial _{x}w\left( t,x\right) \right|
<+\infty .
\end{equation*}
\end{Definition}

\begin{Definition}
\label{df:mildHJBappl1}A function $u\in \Sigma ^{1,\frac{1}{2}}\left( \left[
0,T\right] \times \bar {\mathcal{O}}\right) $ that satisfies the integral
equation (\ref{eq:HJBmildappl1}) for every $\left( t,x\right) \in \left[ 0,T%
\right] \times \bar {\mathcal{O}}$ is called a \emph{mild solution} of the
HJB equation (\ref{eq:HJBappl1}).
\end{Definition}

\begin{Theorem}
\label{th:esolmildappl1} Assume that $B\in C^{0}\left(
\bar{\mathcal{O}}\right) $, $F_{0}\in C^{0}\left( \left[ 0,T\right]
\times \bar{\mathcal{O}}\right) $, $ \phi \in C^{0}\left(
\bar{\mathcal{O}}\right)$ and $\psi=0$. Moreover let $H^{0}$ be
continuous in $\left[ 0,T\right] \times
\bar{\mathcal{O}}\times \mathbb{R}^{n}$ and satisfies (\ref{eq:lipHamappl1}%
). Then there exists a unique mild solution $u\in \Sigma ^{1,\frac{1}{2}%
}\left( \left[ 0,T\right] \times \bar{\mathcal{O}}\right) $of the HJB
equation (\ref{eq:HJBappl1}).
\end{Theorem}

\proof%
%
The proof is a standard application of the contraction mapping principle,
see \cite[Proposition 7.3.4]{Lunardibook} and also \cite{CDPSiamHJ2,CPDE}%
.\hfill \textrm{%
\hbox{\hskip 6pt\vrule width6pt height7pt
depth1pt  \hskip1pt}\bigskip }

Now we want to show that such a mild solution is a strong solution in the
sense that it can be seen as the limit of smooth solutions.

We rewrite here the definition of strong solution for our case since it is
slightly weaker than the one used in Section \ref{PROOFMODEL} due to the
possible singularity of the spatial gradient at $t=0$, see Remark \ref
{rm:singat0proofmodel}.

\begin{Definition}\label{df:strongHJBappl1}
Let $\phi \in C^{0}\left( \bar{\mathcal{O}}\right)$ and $\psi \in
C^{0}\left([0,T]\times \partial{\mathcal{O}}\right)$.
A function $u\in \Sigma ^{1,\frac{1}{2}}\left( %
\left[ 0,T\right] \times \bar{\mathcal{O}}\right) $ is a strong solution of
the equation (\ref{eq:HJBappl1}) if it satisfies the boundary and final
conditions
\begin{equation*}
\begin{array}{l}
u\left( T,x\right) =\phi \left( x\right) ,\qquad x\in \mathcal{O}, \\
u\left( t,x\right) =\psi \left( t,x\right) ,\quad \quad t\in \left[ 0,T%
\right] ,s\in \partial \mathcal{O},
\end{array}
\end{equation*}
and if, for every $\varepsilon >0$, it is a strong solution (in the sense of
Definition \ref{df:solstrong}) of the linear parabolic problem
\begin{eqnarray*}
w_{t}+Lw &=&h,\qquad t\in \left[ 0,T\right] \text{,\qquad }x\in \mathcal{O},
\\
w\left( T-\varepsilon ,x\right)  &=&u\left( T-\varepsilon ,x\right) ,\qquad
x\in \mathcal{O},
\end{eqnarray*}
where $h=\left[ \left\langle F_{0}\left( \cdot ,\cdot \right) ,\partial
_{x}u\left( \cdot ,\cdot \right) \right\rangle +H^{0}\left( \cdot ,\cdot
,\partial _{x}u\left( \cdot ,\cdot \right) \right) \right] $, i.e. if, for
every $\varepsilon >0$ there exist sequences $u_{n}^{\varepsilon }$, $%
h_{n}^{\varepsilon }$, $\phi _{n}^{\varepsilon }$ such that

\begin{enumerate}
\item  for every $n$ $u_{n}^{\varepsilon }$ is a strict solution (i.e. it
belong to $C^{1,2}\left( \left[ 0,T-\varepsilon \right] \times \bar{\mathcal{%
O}}\right) $ and satisfy the equalities on $\left[ 0,T-\varepsilon \right]
\times \bar{\mathcal{O}}$) of the approximating problem
\begin{equation*}
v_{t}+Lv=h_{n}^{\varepsilon },\quad v\left( T-\varepsilon ,x\right) =\phi
_{n}^{\varepsilon }\left( x\right) ;
\end{equation*}

\item  $u_{n}^{\varepsilon }$ converges to $u$ uniformly in $\left[
0,T-\varepsilon \right] \times \bar {\mathcal{O}}$, as $n \to
+\infty$;

\item  $h_{n}^{\varepsilon }$ converges to $h$ uniformly in $\left[
0,T-\varepsilon \right] \times \bar{\mathcal{O}}$, as $n \to
+\infty$;

\item  $\phi _{n}^{\varepsilon } - u\left( T-\varepsilon ,\cdot
\right)  $ converges to zero uniformly in $\bar{\mathcal{O}}$, as $n
\to +\infty$.
\end{enumerate}
\end{Definition}

Given the above Definition \ref{df:strongHJBappl1} we can apply directly
Proposition 4.1.8 (see also Theorem 5.1.11) of \cite{Lunardibook} to get the
following.

\begin{Theorem}
Under the same assumptions of Theorem \ref{th:esolmildappl1} the
mild solution of (\ref{eq:HJBappl1}) is also strong.
\end{Theorem}

\proof%
%
It is enough to apply Proposition 4.1.8 and Theorem 5.1.11 of \cite
{Lunardibook} to the nonhomogenous linear parabolic problem
\begin{equation*}
v_{t}+Lv=h,\quad v\left( T,x\right) =u\left( T-\varepsilon ,x\right) ,
\end{equation*}
where we set
\begin{equation*}
h\left( t,x\right) =\left\langle F_{0}\left( t,x\right) ,\partial
_{x}u\left( t,x\right) \right\rangle +H^{0}\left( t,x,\partial _{x}u\left(
t,x\right) \right) .
\end{equation*}
In fact for every $\varepsilon >0$ such function $h$ belongs to the space $%
C^{0}\left( \left[ 0,T-\varepsilon \right] ;C^{0}\left( \bar{\mathcal{O}}%
\right) \right) =C^{0}\left( \left[ 0,T-\varepsilon \right] \times \bar{%
\mathcal{O}}\right) $.\hfill \textrm{%
\hbox{\hskip 6pt\vrule width6pt height7pt
depth1pt  \hskip1pt}\bigskip }

\begin{Remark}
In the case when $\psi $ is not $0$ the above results still hold
true using the same techniques shown in \cite{Lunardibook}, Theorem
5.1.16, and 5.1.17. We do not do it here for simplicity of
exposition.\hfill\qed
\end{Remark}

\begin{Remark}
Suppose that the operator $L$ can be written in divergence form and that the coefficients
are only Borel measurable and not necessarily continuous; suppose moreover
that the diffusion coefficients are  lower and upper bounded by a constant.
Then, the semigroup has a density with
respect to the Lebesgue measure and it fulfills the classical Aronson
estimates,  see for instance \cite{Aronson, stroock}.

 Fukushima - Dirichlet decomposition for mild or weak solutions to
equations of type (\ref{eq:HJBappl1}) were treated by \cite{BPS,
lejay, rozkosz}. Using such kind of results or possible
generalizations in the spirit of \cite{FRFGpart1}, one could
establish verification theorems related to optimal control problems
even in that framework. \hfill\qed
\end{Remark}

\subsection{The verification theorem\label{FINHORVER}}

Now we prove a verification theorem for the optimal control problem above (%
\ref{eq:SE1}) - (\ref{eq:CFappl1}). The proof is a modification of the proof
given in Section \ref{PROOFMODEL}. The main differences are due to

\begin{enumerate}
\item  a singularity of the first derivative and so a different definition
of strong solution;

\item  the constraint on the set $\mathcal{O}$ and so the presence of
boundary data in $x$ (exit time).
\end{enumerate}

The statement of the verification theorem in this case is the following

\begin{Theorem}
\label{th:verifappl1} Assume that Hypotheses \ref{hp:appl1driftbounded}, \ref
{hp:appl1nondeg}, \ref{hp:appl1control}, \ref{hp:appl1functional}, hold. Let
$v\in \Sigma ^{1,\frac{1}{2}}\left( \left[ 0,T\right] \times \bar {\mathcal{O%
}}\right) $ be a strong solution of (\ref{eq:HJBappl1}). Then

\begin{enumerate}
\item[(i)]  $v\leq V$ on $\left[ 0,T\right] \times \bar {\mathcal{O}}$.

\item[(ii)]  If $z$ is an admissible control at $\left( t,x\right) \in \left[
0,T\right] \times \mathcal{O}$ that satisfies (setting $y\left( s\right)
=y\left( s;t,x,z\right) $), $\mathbb{P}$-a.s.,
\begin{equation*}
H^{0}\left( s,y\left( s\right) ,\partial _{x}v\left( s,y\left( s\right)
\right) \right) =H_{CV}^{0}\left( s,y\left( s\right) ,\partial _{x}v\left(
s,y\left( s\right) \right) ;z\left( s\right) \right) ,
\end{equation*}
for a.e. $s\in \left[ t,T\wedge \tau _{\mathcal{O}}\right] $, then $z$ is
optimal at $\left( t,x\right) $ and $v\left( t,x\right) =V\left( t,x\right) .
$
\end{enumerate}
\end{Theorem}

The proof of this theorem follows as usual by the following fundamental
identity that we state as a lemma.

\begin{Lemma}
Assume that Hypotheses \ref{hp:appl1driftbounded}, \ref{hp:appl1nondeg}, \ref
{hp:appl1control}, \ref{hp:appl1functional}, hold. Let $v\in \Sigma ^{1,%
\frac{1}{2}}\left( \left[ 0,T\right] \times \bar{\mathcal{O}}\right) $ be a
strong solution of (\ref{eq:HJBappl1}). Then, for every $\left( t,x\right)
\in \left[ 0,T\right] \times \bar{\mathcal{O}}$ and $z\in \mathcal{Z}%
_{ad}\left( t\right) $ the following identity holds
\begin{equation}
J\left( t,x;z\right) =v\left( t,x\right)   \notag  \label{eq:fundidappl1}
\end{equation}
\qquad
\begin{equation}
+\mathbb{E}\int_{t}^{\tau _{\mathcal{O}}\wedge T}\left[ -H^{0}\left(
s,y(s),\partial _{x}v\left( s,y\left( s\right) \right) \right)
+H_{CV}^{0}\left( s,y(s),\partial _{x}v\left( s,y\left( s\right) \right)
;z(s)\right) \right] ds,  \label{eq:idfondappl1}
\end{equation}
where $y(s)\overset{def}{=}y\left( s;t,x,z\right) $ is the solution of (\ref
{eq:SE1}) associated with the control $z$.
\end{Lemma}

\proof%
%
Since $v\in \Sigma ^{1,\frac{1}{2}}\left( \left[ 0,T\right] \times \bar{%
\mathcal{O}}\right) $ then for every $\varepsilon >0$ we have $v\in
C^{0,1}\left( \left[ 0,T-\varepsilon \right] \times \bar{\mathcal{O}}\right)
$ and $v$ is a strong solution of the problem
\begin{eqnarray*}
w_{t}+Lw &=&h,\qquad t\in \left[ 0,T\right] \text{,\qquad }x\in \mathcal{O},
\\
w\left( T-\varepsilon ,x\right) &=&v\left( T-\varepsilon ,x\right) ,\qquad
x\in \mathcal{O},
\end{eqnarray*}
where $h=\left[ \left\langle F_{0}\left( \cdot ,\cdot \right) ,\partial
_{x}v\left( \cdot ,\cdot \right) \right\rangle +H\left( \cdot ,\cdot
,\partial _{x}v\left( \cdot ,\cdot \right) \right) \right] $, (see
Definition \ref{df:strongHJBappl1}).

Applying now Theorem \ref{th:Itostrong} and Remark 4.9 of \cite{FRFGpart1}
for the operator $\mathcal{L}_{0}=\partial _{t}+L$ (so now $b=0$, $\sigma =B$
and $b_{1}=F_{0}+F_{1}$), we get that
\begin{equation*}
v\left( \tau _{\mathcal{O}}\wedge \left( T-\varepsilon \right) ,y\left( \tau
_{\mathcal{O}}\wedge \left( T-\varepsilon \right) \right) \right) =v\left(
t,y\left( t\right) \right)
\end{equation*}
\begin{eqnarray*}
&&+\int_{t}^{\tau _{\mathcal{O}}\wedge \left( T-\varepsilon \right)
}\left\langle \partial _{x}v\left( s,y\left( s\right) \right) ,B\left(
y\left( s\right) \right) dW\left( s\right) \right\rangle  \\
&&+\int_{t}^{\tau _{\mathcal{O}}\wedge \left( T-\varepsilon \right)
}H^{0}\left( s,y\left( s\right) ,\partial _{x}v\left( s,y\left( s\right)
\right) \right) ds \\
&&+\int_{t}^{\tau _{\mathcal{O}}\wedge \left( T-\varepsilon \right)
}\left\langle F_{1}\left( s,y\left( s\right) ,z\left( s\right) \right)
,\partial _{x}v\left( s,y\left( s\right) \right) \right\rangle ds.
\end{eqnarray*}
Adding and subtracting $\int_{t}^{\tau _{\mathcal{O}}\wedge \left(
T-\varepsilon \right) }l\left( s,y\left( s\right) ,z\left( s\right) \right)
ds$ (which is always finite since $l$ is bounded), using that $y\left(
t\right) =x$ and taking the expectation we get
\begin{equation*}
\mathbb{E}\left[ \int_{t}^{\tau _{\mathcal{O}}\wedge \left( T-\varepsilon
\right) }l\left( s,y\left( s\right) ,z\left( s\right) \right) ds+v\left(
\tau _{\mathcal{O}}\wedge \left( T-\varepsilon \right) ,y\left( \tau _{%
\mathcal{O}}\wedge \left( T-\varepsilon \right) \right) \right) \right]
=v\left( t,x\right)
\end{equation*}
\begin{equation*}
+\mathbb{E}\int_{t}^{\tau _{\mathcal{O}}\wedge \left( T-\varepsilon \right) }%
\left[ -H^{0}\left( s,y\left( s\right) ,\partial _{x}v\left( s,y\left(
s\right) \right) \right) +H_{CV}^{0}\left( s,y\left( s\right) ,\partial
_{x}v\left( s,y\left( s\right) \right) ;z\left( s\right) \right) \right] ds.
\end{equation*}
Let now $\varepsilon \rightarrow 0^{+}$. By continuity of $v$ (and since by
definition the strong solution satisfies the boundary and final condition)
we get,
\begin{equation*}
\mathbb{E}v\left( \tau _{\mathcal{O}}\wedge \left( T-\varepsilon \right)
,y\left( \tau _{\mathcal{O}}\wedge \left( T-\varepsilon \right) \right)
\right) \rightarrow \mathbb{E}\left[ I_{\left\{ \tau _{\mathcal{O}}\geq
T\right\} }\phi \left( y\left( T\right) \right) +I_{\left\{ \tau _{\mathcal{O%
}}<T\right\} }\psi \left( \tau _{\mathcal{O}},y\left( \tau _{\mathcal{O}%
};t,x,z\right) \right) \right] .
\end{equation*}
Moreover, by the boundedness of $l$ we get
\begin{equation*}
\mathbb{E}\int_{t}^{\tau _{\mathcal{O}}\wedge \left( T-\varepsilon \right)
}l\left( s,y\left( s\right) ,z\left( s\right) \right) ds\rightarrow \mathbb{E%
}\int_{t}^{\tau _{\mathcal{O}}\wedge T}l\left( s,y\left( s\right) ,z\left(
s\right) \right) ds.
\end{equation*}
Finally we recall that, by the definition of the Hamiltonian (\ref
{eq:appl1Hamiltonian}) and thanks to the fact that $v\in \Sigma
^{1,\frac{1}{2}}\left( \left[ 0,T\right] \times
\bar{\mathcal{O}}\right) $,
 for a suitable $C_{0}>0$ we have for every $s\in %
\left[ 0,T\right) $
\begin{eqnarray*}
&&\mathbb{E}\left[ -H^{0}\left( s,y\left( s\right) ,\partial _{x}v\left(
s,y\left( s\right) \right) \right) +H_{CV}^{0}\left( s,y\left( s\right)
,\partial _{x}v\left( s,y\left( s\right) \right) ;z\left( s\right) \right) %
\right] ds \\
&\leq &C_{0}\left( 1+\mathbb{E}\left| \partial _{x}v\left( s,y\left(
s\right) \right) \right| \right) \leq C_{1}\left( T-s\right) ^{-\frac{1}{2}},
\end{eqnarray*}
which allows us to apply the dominated convergence theorem obtaining
\begin{eqnarray*}
&&\mathbb{E}\int_{t}^{\tau _{\mathcal{O}}\wedge \left( T-\varepsilon \right)
}\left[ -H^{0}\left( s,y\left( s\right) ,\partial _{x}v\left( s,y\left(
s\right) \right) \right) +H_{CV}^{0}\left( s,y\left( s\right) ,\partial
_{x}v\left( s,y\left( s\right) \right) ;z\left( s\right) \right) \right] ds
\\
&\rightarrow &\mathbb{E}\int_{t}^{\tau _{\mathcal{O}}\wedge T}\left[
-H^{0}\left( s,y\left( s\right) ,\partial _{x}v\left( s,y\left( s\right)
\right) \right) +H_{CV}^{0}\left( s,y\left( s\right) ,\partial _{x}v\left(
s,y\left( s\right) \right) ;z\left( s\right) \right) \right] ds.
\end{eqnarray*}
All the above yields then the claim (\ref{eq:idfondappl1}).\hfill
\hbox{\hskip 6pt\vrule width6pt
height7pt depth1pt  \hskip1pt}\bigskip

\textbf{Proof of Theorem \ref{th:verifappl1}. } By the definition of $H^{0}$
and $H_{CV}^{0}$, for every $z\in \mathcal{Z}_{ad}\left( t\right) $, $x\in
\mathcal{O}$ the following inequality holds $\mathbb{P}$-almost surely
(setting $y\left( s\right) =y\left( s;t,x,z\right) $)
\begin{equation}
H^{0}\left( s,y\left( s\right) ,\partial _{x}v\left( s,y\left( s\right)
\right) \right) -H_{CV}^{0}\left( s,y\left( s\right) ,\partial _{x}v\left(
s,y\left( s\right) \right) ;z\left( s\right) \right) \geq 0,\quad \text{for }%
a.e.s\in \left[ t,T\wedge \tau _{\mathcal{O}}\right] .
\label{eq:distrivappl1}
\end{equation}
The rest of the proof follows exactly as the one of Theorem \ref
{th:VTintronostro}. \hfill
\hbox{\hskip 6pt\vrule width6pt
height7pt depth1pt  \hskip1pt}\bigskip

\begin{Remark}
Arguing as in Subsection \ref{FEEDBACK} it is possible to prove a necessary
condition and the existence of weakly optimal feedbacks.\hfill
\hbox{\hskip 6pt\vrule width6pt height7pt
depth1pt  \hskip1pt}\bigskip
\end{Remark}

\section{Application 2: a class of problems with degenerate diffusion\label%
{APPL2DEG}}

Now we consider our model problem focusing the case where the HJB equation
admits solutions that are $C^{0,1}$ but not $C^{1,2}$. This may occur for
instance in the case when the matrix $BB^{\ast }$ is degenerate in the sense
that there exists at least one point $x\in \mathbb{R}^{n}$ such that $%
BB^{\ast }$ is not invertible.

Indeed such degeneracy means that the HJB equations (\ref{eq:HJBristretta})
and (\ref{eq:HJBell}) are degenerate (parabolic or elliptic). For such
equations a general regularity theory like the one available in the
uniformly parabolic or elliptic case is not known. So it is difficult to
obtain existence of $C^{1,2}$ ($C^{2}$) solutions in this case. There is
more hope to obtain solutions with regularity $C^{0,1}$ ($C^{1}$) or at
least continuous with existence of the space derivatives in the weak sense
(this case is not covered in the present paper but is the subject of our
current research). This means that, on one hand, this kind of problems is a
source of possible sharp applications of our verification Theorems \ref
{th:VTintronostro} and \ref{th:VTinfhor}. On the other hand in this case it
is difficult to have Hypothesis \ref{hp:exsolHJBvera}-(ii) to be satisfied.
So one needs to have convergence of the derivatives in the sense of
Hypothesis \ref{hp:exsolHJBvera}-(i), which may be not easy to check.\ We
give an example when this can be done in the last part of this section.

\begin{Remark}
\label{rm:solstrongviasemigroup}A possible methodology to prove existence of
$C^{0,1}$ solutions of the parabolic degenerate equation (\ref
{eq:HJBristretta}) is the following (consider the case where all data are
autonomous). Take the degenerate elliptic operator
\begin{equation*}
L_{1}\eta \left( x\right) =\left\langle F_{0}\left( x\right) ,\partial
_{x}\eta \left( x\right) \right\rangle +\frac{1}{2}\text{Tr }\left[ B^{\ast
}\left( x\right) \partial _{xx}\eta \left( x\right) B\left( x\right) \right]
,
\end{equation*}
defined in $C^{2}\left( \mathbb{R}^{n}\right) $. We can write the HJB
equation (\ref{eq:HJBristretta}) as
\begin{equation*}
v_{t}\left( t,x\right) +L_{1}v\left( t,x\right) +H^{0}\left( x,\partial
_{x}v\left( t,x\right) \right) =0,
\end{equation*}
\begin{equation*}
v\left( T,x\right) =\phi \left( x\right) .
\end{equation*}
At this point, if the operator $L_{1}$ is ``sufficiently good'', e.g. if it
generates a smoothing semigroup (see the papers \cite{GV,GMV} on this) one
can try to apply some perturbation arguments finding existence of a strong
solutions. Also Remark \ref{rm:singat0proofmodel} applies.\hfill
\hbox{\hskip 6pt\vrule
width6pt height7pt depth1pt  \hskip1pt}\bigskip
\end{Remark}

To illustrate the above situation we present a one dimensional example
coming from optimal advertising models. For this reason the problem is taken
with the $\sup $, so also Hamiltonians are different from the rest of the
paper. The proofs are only sketched for brevity.

The state variable $y$ is the goodwill (i.e. the number of people that are
aware of the product), the control variable $z$ is the investment in
advertising and one maximizes the profit from selling the product on a given
finite horizon; see on this e.g. \cite{Grosset,Marinelli} and the references
therein.

The state space is $\mathbb{R}$, the control space is $U=\mathbb{R}^{+}$.
The noise space is $\mathbb{R}$ and $W$ is a standard 1-dimensional Brownian
motion. Given an initial point $x\in \mathbb{R}$ and parameters $\alpha
,\beta >0$, the state equation is
\begin{equation}
\left\{
\begin{array}{c}
dy\left( s\right) =\left[ -\alpha y\left( s\right) +z\left( s\right) \right]
ds+\beta y\left( s\right) dW\left( s\right)  \\
y\left( t\right) =x.\qquad \qquad \qquad \qquad \qquad \qquad \qquad
\end{array}
\right.   \label{eq:SE1dim1}
\end{equation}
Given $\rho >0$, $c:\mathbb{R}^{+}\mathbb{\rightarrow R}^{+}$, continuous,
increasing and convex and such that
\begin{equation}
\lim_{z\rightarrow +\infty }\frac{c\left( z\right) }{z}=+\infty ,
\label{eq:coercivec}
\end{equation}
$h:\mathbb{R}\rightarrow \mathbb{R}$, continuous and strictly increasing, we
maximize the profit functional
\begin{equation*}
J(t,x;z)=\mathbb{E}\left\{ \int_{t}^{T}-c\left( z\left( s\right) \right)
ds+h\left( y\left( T\right) \right) \right\} ,
\end{equation*}
over all controls $z\in \mathcal{Z}_{ad}\left( t\right) $ where
\begin{equation*}
\mathcal{Z}_{ad}(t)=\{z:\mathcal{T}_{t}\times \Omega \rightarrow U,\quad
\text{measurable, a.s. locally integrable, }(\mathcal{F}_{s})\text{-adapted}%
\}.
\end{equation*}
Here $y  =y\left( \cdot ;x,z\right) =y\left( \cdot
;t,x,z\right) $ is the strong solution of the equation (\ref{eq:SE1dim1}) on
$\left[ 0,+\infty \right) $. Here Hypotheses \ref{hp:introstateeq}, \ref
{hp:introcost} and \ref{hp:Hamiltonian1} are clearly satisfied. The value
function is
\begin{equation*}
V(t,x)=\sup \left\{ J(t,x;z):\;z\in \mathcal{Z}_{ad}\left( t\right) \right\}
\end{equation*}
and the corresponding Hamilton-Jacobi-Bellman equation (HJB from now on)
reads as follows
\begin{equation}
\left\{
\begin{array}{l}
\partial _{t}v\left( t,x\right) +\frac{1}{2}\mathrm{\beta }^{2}x^{2}\partial
_{xx}v(t,x)-\left\langle \alpha x,\partial _{x}v\left( t,x\right)
\right\rangle +H^{0}\left( \partial _{x}v\left( t,x\right) \right) =0, \\
\qquad \qquad \qquad \qquad \qquad \qquad \qquad \qquad \qquad \qquad \qquad
t\in \left[ 0,T\right] \text{,\qquad }x\in \mathbb{R}, \\
v\left( T,x\right) =h\left( x\right) ,\qquad x\in \mathbb{R,}
\end{array}
\right.   \label{eq:HJBes1appl2}
\end{equation}
where the Hamiltonian $H^{0}$ is given by
\begin{equation*}
H^{0}\left( p\right) =\sup_{z\geq 0}\left\{ zp-c\left( z\right) \right\}
=c^{\ast }\left( p\right) ,
\end{equation*}
(where $c^{\ast }$ is the Legendre transform of $c$) and it is always finite
and continuous thanks to (\ref{eq:coercivec}).

\begin{Remark}
\label{rm:betterthanalternative}We stress the fact that, via the approach
described in Remark \ref{rm:solstrongviasemigroup}, we can prove, eventually
restricting the assumptions on data, that there exists a strong solution of
the HJB equation (\ref{eq:HJBes1appl2}). In fact, the operator
\begin{equation*}
L_{1}\eta \left( x\right) =\left\langle \alpha x,\partial _{x}\eta \left(
x\right) \right\rangle +\frac{1}{2}\beta ^{2}x^{2}\partial _{xx}\eta \left(
x\right) ,
\end{equation*}
generates an analytic semigroup (see e.g. \cite{GMV} and the references
therein) in the space $C^{0}$ with domain $C^{2}$ with suitable weights.
This allows us, using a perturbation method similar to the one of \cite{JMAA}%
, to show, under some restriction on the data, the existence of a mild
solution of HJB equation (\ref{eq:HJBes1appl2}), i.e. a solution of the
integral equation
\begin{equation*}
v\left( t,\cdot \right) =e^{\left( T-t\right) L_{1}}h+\int_{t}^{T}e^{\left(
T-s\right) L_{1}}H^{0}\left( \partial _{x}v\left( s,\cdot \right) \right) ds
\end{equation*}
belonging to the space $C^{0,1}$ with suitable weights. Such mild solution
turns out to be strong in the sense of Definition \ref{df:solHJBstrong}
easily: by simply taking suitable approximations of the initial datum $h$
and of $H^{0}\left( \partial _{x}v\left( s,\cdot \right) \right) $, as it is
done again in \cite{JMAA} or in \cite{Lunardibook}, one finds solutions $%
u_{n}$ of approximating problems that converges to $u$. However it is not
trivial to see if for such approximating sequence we also have the
convergence required in Hypothesis \ref{hp:exsolHJBvera}-(i). Below we see
that it is true in a simple case.
\hbox{\hskip 6pt\vrule width6pt height7pt depth1pt
\hskip1pt}\bigskip
\end{Remark}

Choosing data in a suitable way we can write explicitely a strong solution
which satisfies our assumptions but is not $C^{2}$ in space. Let $c\left(
z\right) =z^{1+\eta }$ with $\eta \in \left( 0,1\right) $, and $h\left(
x\right) =\left| x\right| ^{1+\eta }sgn\left( x\right) $. Then the
Hamiltonian is
\begin{equation*}
H^{0}\left( p\right) =\eta \left( \frac{\left[ p\right] ^{+}}{1+\eta }%
\right) ^{1+\frac{1}{\eta }},
\end{equation*}
with the maximum point
\begin{equation}
z=\left( \frac{\left[ p\right] ^{+}}{1+\eta }\right) ^{\frac{1}{\eta }}.
\label{eq:puntomaxappl2es1}
\end{equation}

If we look for a solution of the form $w\left( t,x\right) =f\left( t\right)
\cdot \left| x\right| ^{1+\eta }$ we see that the function
\begin{equation*}
v\left( t,x\right) =\left\{
\begin{array}{cc}
a\left( t\right) \left| x\right| ^{1+\eta }, & \text{for }x>0, \\
0, & \text{for }x=0, \\
b\left( t\right) \left| x\right| ^{1+\eta }, & \text{for }x<0,
\end{array}
\right.
\end{equation*}
satisfies for every $x\neq 0$ the HJB equation (\ref{eq:HJBes1appl2}) above
if $a>0$ and $b<0$ are, respectively, the solution of the Cauchy problem
\begin{equation*}
a^{\prime }\left( t\right) =-\left[ \frac{1}{2}\beta ^{2}\eta \left( 1+\eta
\right) -\alpha \left( 1+\eta \right) \right] a\left( t\right) -\eta a\left(
t\right) ^{1+\frac{1}{\eta }},\qquad a\left( T\right) =1,
\end{equation*}
and, respectively,
\begin{equation*}
b^{\prime }\left( t\right) =\left[ \frac{1}{2}\beta ^{2}\eta \left( 1+\eta
\right) -\alpha \left( 1+\eta \right) \right] b\left( t\right) ,\qquad
b\left( T\right) =-1.
\end{equation*}
Now the second Cauchy problem admits a unique global solution on $\left[ 0,T%
\right] $, which is always strictly negative. The first does the same under
suitable conditions on the data (e.g. that $\frac{1}{2}\beta ^{2}\eta
<\alpha $): assume from now on that this is the case. We can easily see that
$v\in C^{0,1}\left( \left[ 0,T\right] \times \mathbb{R}\right) $. Moreover
we may prove, by approximating $v$ by $v_{n}=v\cdot \theta _{n}$ (where $%
\theta _{n}$ is a suitable cut-off function such that $\theta _{n}=0$ for $%
\left| x\right| \leq \frac{1}{n}$ and $\theta _{n}=1$ for $\left| x\right|
\geq \frac{2}{n}$) that $v$ is a strong solution of (\ref{eq:HJBes1appl2})
and that $\partial _{x}v_{n}\rightarrow \partial _{x}v$ uniformly on compact
sets. This fact says that we have a strong solution to which standard
verification do not apply but which falls into our assumptions. Indeed since
the convergence of the derivative is uniform on compact sets also another
technique can be applied (see Subsection \ref{COMPAPPR}).

Remark that not only Theorem \ref{th:VTintronostro} but also Proposition \ref
{pr:feedbackvero} applies. Indeed the maximum of the Hamiltonian is reached
in the unique point $z$ given by (\ref{eq:puntomaxappl2es1}), so we can set,
with the notation of Proposition \ref{pr:feedbackvero}
\begin{equation*}
G_{0}(t,x,p)=G_{0}(p)=\left( \frac{\left[ p\right] ^{+}}{1+\eta }\right) ^{%
\frac{1}{\eta }}.
\end{equation*}
Recalling that
\begin{equation*}
\partial _{x}v\left( t,x\right) =\left\{
\begin{array}{cc}
a\left( t\right) \left( 1+\eta \right) \left| x\right| ^{\eta }, & \text{for
}x>0, \\
0, & \text{for }x=0, \\
-b\left( t\right) \left( 1+\eta \right) \left| x\right| ^{\eta }, & \text{%
for }x<0,
\end{array}
\right.
\end{equation*}
we have
\begin{equation*}
G(t,x)=G_{0}(\partial _{x}v\left( t,x\right) )=\left\{
\begin{array}{cc}
a\left( t\right) ^{\frac{1}{\eta }}\left| x\right|  & \text{for }x>0 \\
0 & \text{for }x=0 \\
\left[ -b\left( t\right) \right] ^{\frac{1}{\eta }}\left| x\right|  & \text{%
for }x<0
\end{array}
\right.
\end{equation*}
It is easy to check that $G$ is an admissible feedback map (along with
Definition \ref{df:feedback}). In fact, setting, for $x>0$%
\begin{equation*}
z\left( s\right) =a\left( s\right) ^{\frac{1}{\eta }}\left| y\left( s\right)
\right| ,
\end{equation*}
the closed loop equation is
\begin{equation*}
\left\{
\begin{array}{c}
dy\left( s\right) =\left[ -\alpha y\left( s\right) +a\left( s\right) ^{\frac{%
1}{\eta }}\left| y\left( s\right) \right| \right] dt+\beta y\left( s\right)
dW\left( s\right)  \\
y\left( t\right) =x.\qquad \qquad \qquad \qquad \qquad \qquad \qquad
\end{array}
\right.
\end{equation*}
This equation always features existence and uniqueness of a strong solution $%
y_{G,t,x}$ which is a.s. strictly positive and so $z_{G,t,x}$ is admissible
and $z_{G,t,x}\left( s\right) =G(s,y_{G,t,x}\left( s\right) )$ for $s\in
\left[ t,T\right] $. The same argument can be applied to the cases $x<0$ and
$x=0$. This means that, from Proposition \ref{pr:feedbackvero}, for every
fixed $(t,x)\in \lbrack 0,T]\times \mathbb{R}$ the couple $%
(z_{G,t,x},y_{G,t,x})$ (again with the notations of Proposition \ref
{pr:feedbackvero}) is optimal. Moreover the optimal control (state) is
unique. In particular, when $x>0$ ($x<0$) then we have a unique optimal
couple where the optimal state is a.s. positive (negative). When $x=0$ the
zero strategy is optimal. These facts cannot be immediately deduced using
other verification theorems. Of course one could argue with ad hoc arguments
applied to this problem but this is another story. We observe that this also
implies that $v$ is the value function.

\begin{Remark}
The procedure outlined in this example may be applied to other kind
of problems with bilinear state equation and current cost with
growth $1+\eta $ ($\eta \in \left( 0,1\right) $), eventually with
infinite horizon. We think e.g. to optimal investment models where
adjustment cost are not quadratic but of growth $1+\eta $ (see e.g.
\cite{Takayama}). Also one could treat similar models when no
explicit solution can be found but it is possible to prove that
there exists a strong $C^{0,1}$ solution of the HJB equation.

Finally we mention two possible extensions of the problem introduced above
that are interesting for economic applications (namely growth theory,
optimal investment, optimal portfolio models): the case when state
constraints arise and the case when the space gradient of the value function
presents some singularities. Take for example the advertising example with $%
\eta <0$. These cases do not fall in our assumptions but we think that our
procedure can be extended to cover them in a future work.\hfill
\hbox{\hskip 6pt\vrule width6pt height7pt
depth1pt  \hskip1pt}\bigskip
\end{Remark}

\bigskip

\section{Comparison with other verification results\label{COMPARISON}}

In this section we compare our results with other known verification
techniques for non smooth solution of the HJB equation. A discussion about
this matter has already been done in Section \ref{INTROVER}. Here we give a
more detailed analysis.

\subsection{The classical result}

We start by recalling that the classical verification theorems for
stochastic optimal control problems (see e.g. \cite[p.163]{FS}) adapted to
our case are perfectly equal to our Theorems \ref{th:VTintronostro}, \ref
{th:VTinfhor} except for the following facts:

\begin{itemize}
\item  they assume, in the finite horizon case, $v\in C^{1,2}\left( \left[
0,T\right] \times \mathbb{R}^{n}\right) $ instead of $v\in C^{0,1}\left( %
\left[ 0,T\right] \times \mathbb{R}^{n}\right) $;

\item  they assume, in the infinite horizon case, $v\in C^{2}\left( \mathbb{R%
}^{n}\right) $ instead of $v\in C^{1}\left( \mathbb{R}^{n}\right) $.
\end{itemize}

In Section \ref{APPL1EXIT} we have seen an example where $v$ is surely $%
C^{0,1}\left( \left[ 0,T\right] \times \mathbb{R}^{n}\right) $ but it is not
clear if it is also $C^{1,2}\left( \left[ 0,T\right] \times \mathbb{R}%
^{n}\right) $. Moreover, in Section \ref{APPL2DEG} we have seen a one
dimensional example where $v$ is surely $C^{0,1}\left( \left[ 0,T\right]
\times \mathbb{R}\right) $ and not $C^{1,2}\left( \left[ 0,T\right] \times
\mathbb{R}\right) $.

\subsection{The approximation result\label{COMPAPPR}}

By approximation result we mean the verification theorem, proved e.g. in
\cite{CPDE} in a special infinite dimensional case, that is proved using the
approach that we have briefly recalled in Remark \ref{rm:suffperconvder}.
The statement of such result is completely similar to our Theorems \ref
{th:VTintronostro} or \ref{th:VTinfhor}. The only difference is that to
prove this theorem one needs to know that the function $v$ is a strong
solution of HJB equation in a more restrictive sense (we may call them
``stronger'' solutions) pointed out in Remark \ref{rm:suffperconvder}: it is
required that also the space derivative of the approximating sequence
converges uniformly on compact sets, i.e. $\partial _{x}v_{n}\rightarrow
\partial _{x}v$. In the results of the present paper in place of this we
have a weaker requirement, namely that, either Hypothesis \ref
{hp:exsolHJBvera}-(i), or \ref{hp:exsolHJBvera}-(ii) is satisfied. Regarding
this point we note that:

\begin{itemize}
\item  in Application 1 the fact that the solution is strong follows
directly from theorems on semigroups stated in \cite{Lunardibook}. Here
Hypothesis \ref{hp:exsolHJBvera}-(i) holds by the nondegeneracy of the
diffusion coefficient while the convergence of the space derivative in such
cases is not trivial at all and may be not true.

\item  in Application 2, in the one dimensional example, as pointed out in
Remark \ref{rm:betterthanalternative}, it is straightforward to prove that $%
u $ is a strong solution in the sense of Definition \ref{df:solHJBstrong}.
Here the diffusion coefficient is degenerate so to apply our results we need
Hypothesis \ref{hp:exsolHJBvera}-(ii). This is not trivial to check in
general but is in any case easier than the uniform convergence of the
derivatives. In the special case where we calculate the value function the
uniform convergence holds so also the approach of \cite{CPDE} may be used.
\end{itemize}

We finally point out that, to avoid the requirement of the
convergence of the derivative on compact sets we need to use a
completely different (and much more complex) approach based on the
Fukushima-Dirichlet decomposition and on the representation result
(Theorem \ref{th:Itostrong}) proved in our companion paper
\cite{FRFGpart1}.

\subsection{The viscosity solution result}

When the HJB equation admits a viscosity solution it is possible to prove a
verification theorem which is very general since it deals with only
continuous solutions of HJB, and with the case when also the diffusion
coefficient is controlled, but which is less useful when we know that we
have strong solutions. The theorem, in the model case studied in this paper,
is the following (see e.g. \cite{GSZstochver,LYZ,Zhoubook} for the proof and
precise definition of superdifferentials $D_{t+,x}^{1,2,+}v(s,y^{\ast }(s))$
used here).

\begin{Theorem}
\label{th:vervisc}Consider the problem (\ref{eq:stateintro}) - (\ref
{eq:cfintro}) with the following assumptions. The data $F_{0}$, $F_{1}$, $B$%
, $l$, $\phi $ are all uniformly continuous. Moreover they are Lipschitz
continuous in the variable $x$ uniformly with respect to the other
variables. Moreover, for a suitable constant $M>0$,
\begin{equation*}
\left| F_{0}\left( t,0\right) \right| ,\left| F_{1}\left( t,0,z\right)
\right| ,\left| B\left( t,0\right) \right| ,\left| l\left( t,0,z\right)
\right| ,\left| \phi \left( 0\right) \right| \leq M.
\end{equation*}
Finally the control space $U$ is a Polish space and the set of admissible
controls $\mathcal{Z}_{ad}\left( t\right) $ is \ given by all measurable and
adapted processes on $\left[ t,T\right] $ with values in $U$.

\textit{L}et $v\in C^{0}([0,T]\times \mathbb{R}^{n})$ be a viscosity
solution of the HJB equation (\ref{eq:HJBintropar}). Then:

\begin{itemize}
\item  $v\left( t,x\right) \leq V\left( t,x\right) $, for any $\left(
t,x\right) \in \lbrack 0,T)\times \mathbb{R}^{n}$.

\item  Let\ $\left( y^{\ast }(\cdot ),z^{\ast }(\cdot )\right) $\ be\ a\
given\ admissible\ pair\ for\ the\ problem starting at $\left( t,x\right) $.
Suppose\ that\ there\ exists\ $\left( p^{\ast },q^{\ast },Q^{\ast }\right)
\in L_{\mathcal{F}}^{2}\left( t,T;\mathbb{R}\right) \times L_{\mathcal{F}%
}^{2}\left( t,T;\mathbb{R}^{n}\right) \times L_{\mathcal{F}}^{2}\left(
t,T;S^{n}\right) $\ such\ that\ for\ a.e.$s\in \lbrack t,T]$,\
\begin{equation*}
(p^{\ast }(s),q^{\ast }(s),Q^{\ast }(s))\in D_{t+,x}^{1,2,+}v(s,y^{\ast
}(s)),\;\mathbb{P}-a.s.
\end{equation*}
and
\begin{eqnarray*}
&&p^{\ast }(s)-\frac{1}{2}\mathrm{Tr}\left[ B^{\ast }\left( s,x^{\ast
}(s)\right) Q^{\ast }\left( s\right) B\left( s,x^{\ast }(s)\right) \right]
-\left\langle F_{0}\left( s,x^{\ast }(s)\right) ,q^{\ast }\left( s\right)
\right\rangle  \\
&=&H_{CV}^{0}(s,x^{\ast }(s),q^{\ast }(s);z^{\ast }(s))=H^{0}(s,x^{\ast
}(s),q^{\ast }(s)),\qquad \;\mathbb{P}-a.s.
\end{eqnarray*}
then\ $(y^{\ast }(\cdot ),z^{\ast }(\cdot ))$\ is\ an\ optimal\ pair\ for\
the\ problem starting at $\left( t,x\right) $.
\end{itemize}
\end{Theorem}

If we know that there exists a strong solution $v\in C^{0,1}$ then
to apply the above theorem we need to know that $v$ is also a
viscosity solution. This is not difficult as the concept of
viscosity solution is more general in a wide range of cases (this is
true since classical solutions are always viscosity solutions and
limits of viscosity solutions are still viscosity solutions, see
e.g. \cite{FS}); if so then we know that it must be
\begin{equation*}
q^{\ast }\left( s\right) =\partial _{x}v\left( s,y^{\ast }\left( s\right)
\right) .
\end{equation*}
\textit{The above sufficient condition states then that:}

if there\ exists\ $\left( p^{\ast },q^{\ast },Q^{\ast }\right) \in L_{%
\mathcal{F}}^{2}\left( t,T;\mathbb{R}\right) \times L_{\mathcal{F}%
}^{2}\left( t,T;\mathbb{R}^{n}\right) \times L_{\mathcal{F}}^{2}\left(
t,T;S^{n}\right) $\ such\ that\ for\ a.e.$s\in \lbrack t,T]$,\
\begin{equation*}
(p^{\ast }(s),q^{\ast }\left( s\right) ,Q^{\ast }(s))\in
D_{t+,x}^{1,2,+}v(s,y^{\ast }(s)),\;\mathbb{P}-a.s.
\end{equation*}
(so $q^{\ast }\left( s\right) =\partial _{x}v\left( s,y^{\ast }\left(
s\right) \right) $) and
\begin{eqnarray*}
&&p^{\ast }(s)-\frac{1}{2}\mathrm{Tr}\left[ B^{\ast }\left( s,x^{\ast
}(s)\right) Q^{\ast }\left( s\right) B\left( s,x^{\ast }(s)\right) \right]
-\left\langle F_{0}\left( s,x^{\ast }(s)\right) ,\partial _{x}v\left(
s,y^{\ast }\left( s\right) \right) \right\rangle  \\
&=&H_{CV}^{0}\left( s,x^{\ast }(s),\partial _{x}v\left( s,y^{\ast }\left(
s\right) \right) ;z^{\ast }(s)\right) =H^{0}\left( s,x^{\ast }(s),\partial
_{x}v\left( s,y^{\ast }\left( s\right) \right) \right) ,\;\mathbb{P}-a.s.
\end{eqnarray*}
then\ $(y^{\ast }(\cdot ),z^{\ast }(\cdot ))$\ is\ an\ optimal\ pair\ for\
the\ problem starting at $\left( t,x\right) $.

\smallskip

\textit{Our sufficient conditions instead states the following:}

if for\ a.e. $\!s\in \lbrack t,T]$,
\begin{equation*}
H_{CV}^{0}\left( s,x^{\ast }(s),\partial _{x}v\left( s,y^{\ast
}\left( s\right) \right) ;z^{\ast }(s)\right) =H^{0}\left( s,x^{\ast
}(s),\partial _{x}v\left( s,y^{\ast }\left( s\right) \right) \right)
,\;\mathbb{P}-a.s.,
\end{equation*}
then\ $(y^{\ast }(\cdot ),z^{\ast }(\cdot ))$\ is\ an\ optimal\ pair\ for\
the\ problem starting at $\left( t,x\right) $.

\smallskip

The main difference is that in the Theorem \ref{th:vervisc} above we need to
know about the existence of a couple of processes $p^{\ast }$ and $Q^{\ast }$
(that play the role of time derivative and second space derivative) which
are in fact not needed in our statement. This is not a trivial difficulty in
general. So our statement shows that, when strong solutions exist, one can
get sharper sufficient conditions. Similar and even stronger considerations
can be done for the necessary conditions and the existence of optimal
feedbacks. In particular the analogous of the necessary condition stated by
Proposition \ref{pr:feedback2} in the viscosity setting says only that for
every possible $(p^{\ast }(s),q^{\ast }(s),Q^{\ast }(s))\in
D_{t+,x}^{1,2,+}v(s,y^{\ast }(s)),\;\mathbb{P}-a.s.$ (so $q^{\ast }\left(
s\right) =\partial _{x}v\left( s,y^{\ast }\left( s\right) \right) $) we have
\begin{eqnarray*}
&&p^{\ast }(t)-\frac{1}{2}\mathrm{Tr}\left[ B^{\ast }\left( t,x^{\ast
}(t)\right) Q^{\ast }\left( t\right) B\left( t,x^{\ast }(t)\right) \right]
-\left\langle F_{0}\left( t,x^{\ast }(t)\right) ,\partial _{x}v\left(
t,y^{\ast }\left( t\right) \right) \right\rangle  \\
&\leq &H_{CV}^{0}\left( t,x^{\ast }(t),\partial _{x}v\left( t,y^{\ast
}\left( t\right) \right) ;z^{\ast }(t)\right)
\end{eqnarray*}
and this \textit{do not} imply a priori that
\begin{equation*}
H_{CV}^{0}\left( t,x^{\ast }(t),\partial _{x}v\left( t,y^{\ast }\left(
t\right) \right) ;z^{\ast }(t)\right) =H^{0}\left( t,x^{\ast }(t),\partial
_{x}v\left( t,y^{\ast }\left( t\right) \right) \right) ,\;\mathbb{P}-a.s
\end{equation*}
since it may happen that
\begin{eqnarray*}
&&p^{\ast }(t)-\frac{1}{2}\mathrm{Tr}\left[ B^{\ast }\left( t,x^{\ast
}(t)\right) Q^{\ast }\left( t\right) B\left( t,x^{\ast }(t)\right) \right]
-\left\langle F_{0}\left( t,x^{\ast }(t)\right) ,\partial _{x}v\left(
t,y^{\ast }\left( t\right) \right) \right\rangle  \\
&<&H^{0}\left( t,x^{\ast }(t),\partial _{x}v\left( t,y^{\ast }\left(
t\right) \right) \right)
\end{eqnarray*}
for all possible $(p^{\ast }(s),q^{\ast }(s),Q^{\ast }(s))\in
D_{t+,x}^{1,2,+}v(s,y^{\ast }(s)),\;\mathbb{P}-a.s.$. So, concerning
necessary conditions we may say that the difference between our result and
the viscosity solution result is wider.

\begin{Remark}
We must also note another point. The assumptions made in Theorem \ref
{th:vervisc} about the data are much more restricitve than our assumptions
made for the model problem of Section \ref{APPL1EXIT}. We think that Theorem
\ref{th:vervisc} could be extended to more general situations but this
extension is not trivial and not known at this stage.

Of course the importance of the verification theorem for viscosity solutions
is that it covers the cases when the control enters also in the diffusion
coefficients and when the viscosity solutions has no further regularity. So
it is in a sense natural that such result is more restrictive than ours when
applied to cases when the control enters only in the drift and a strong
solution exists.\hfill
\hbox{\hskip 6pt\vrule width6pt height7pt
depth1pt  \hskip1pt}\bigskip
\end{Remark}

\subsection{The backward equations result}

This is a verification technique introduced for instance in
\cite{Quenez} and \cite{FTHJB} that is used when it is possible to
find a mild solution to the HJB equation (in the sense recalled in
Subsection \ref{STRONGHJB}) and when it can be represented via the
solution of a forward - backward system, see on this
\cite{FTHJB,MaYong}. The most recent and complete result is given in
\cite{FTHJB} in the infinite dimensional case and we take the
setting from it adapting the statement to our problem. Similar
results in other context are given e.g. in \cite{Quenez} (see also
the references therein).

Consider then the problem (\ref{eq:stateintro}) - (\ref{eq:cfintro}) with
the following assumptions.

\begin{enumerate}
\item  The state space is $X=\mathbb{R}^{n}$ while the control space is $U$,
a bounded subset of $\mathbb{R}^{m}$.

\item  The data $F_{0}$, $B$ are continuously differentiable with bounded
first derivative in the variable $x$ uniformly with respect to $t$.

\item  $F_{1}\left( t,x,z\right) =B\left( t,x\right) R\left( t,x\right) z$
where $R:\left[ 0,T\right] \times X\rightarrow L\left( U,X\right) $ is
continuously differentiable in $x$, bounded with its space derivative.

\item  $l$ is continuous and polynomially growing.

\item  $\phi $ is continuously differentiable with bounded first derivatives.

\item  Defining the modified current value Hamiltonian $K_{CV}^{0}\left(
t,x,q;z\right) =l\left( t,x,z\right) +\left\langle q,z\right\rangle _{%
\mathbb{R}^{m}}$ we assume that $K^{0}=\inf K_{CV}^{0}\left( t,x,q;z\right) $
is measurable, continuously differentiable in $\left( x,p\right) $ with
bounded derivative with respect to $p$ and with polynomially growing
derivative with respect to $x$.

\item  There exists a unique minimum point of the modified current value
Hamiltonian $K_{CV}^{0}$ given by $\Gamma \left( t,x,p\right) $ where $%
\Gamma $ is continuous$.$
\end{enumerate}

The result is the following.

\begin{Theorem}
\label{th:FT}Assume Hypotheses 1-7 above. Fix $\left( t,x\right) \in \left[
0,T\right] \times X$. For all admissible control strategies $z$ starting at $%
\left( t,x\right) $ we have $J(t,x;z)\geq v(t,x)$ and the equality holds if
and only if the following feedback law is verified by $z$ and $y\left( \cdot
;t,x,z\right) $:
\begin{equation}
z(s)=\Gamma (s,y\left( s\right) ,R(s,y\left( s\right) )^{\ast }G(s,y\left(
s\right) )^{\ast }\partial _{x}v(s,y\left( s\right) )),\quad \mathbb{P-}%
a.s.\;for\;a.e.\;s\in \lbrack t,T].  \label{leggefeedback}
\end{equation}
Finally there exists at least an admissible control strategy $z$ starting at
$\left( t,x\right) $ for which (\ref{leggefeedback}) holds. In such a system
the closed loop equation:
\begin{equation}
\left\{
\begin{array}{lll}
dy\left( s\right) & = & F_{0}\left( s,y\left( s\right) \right) +B\left(
s,y\left( s\right) \right) R\left( s,y\left( s\right) \right) \Gamma
_{1}\left( s,y\left( s\right) \right) ds \\
&  & +B\left( s,y\left( s\right) \right) dW_{s},\qquad s\in \lbrack t,T], \\
y\left( t\right) & = & x\in X,
\end{array}
\right.  \label{equazioneclosedloop}
\end{equation}
where we set $\Gamma _{1}\left( s,y\left( s\right) \right) =\Gamma \left(
s,y\left( s\right) ,R\left( s,y\left( s\right) \right) ^{\ast }B\left(
s,y\left( s\right) \right) ^{\ast }\nabla _{x}v\left( s,y\left( s\right)
\right) \right) $, admits a solution $\bar{y}$ and if $\bar{z}(s)=\Gamma
\left( s,\bar{y}\left( s\right) ,R(s,\bar{y}\left( s\right) )^{\ast }G(s,%
\bar{y}\left( s\right) )^{\ast }\nabla _{x}v\left( s,\bar{y}\left( s\right)
\right) \right) $ then the couple $(\bar{z},\bar{y})$ is optimal for the
control problem.
\end{Theorem}

We point out the following.

\begin{itemize}
\item  The dependence on the control in the state equation is linear and
through the same operator $B$ that drives the noise, moreover the space $U$
is compact. This means that Hypothesis \ref{hp:exsolHJBvera} (ii) is always
satisfied. So the backward equation result Theorem \ref{th:FT} cannot cover
cases when only part (i) of Hypothesis \ref{hp:exsolHJBvera} is satisfied,
e.g. the one dimensional example of Section \ref{APPL2DEG}. This is not a
casual restriction but it seems to be a limitation due to the technique used
in the proof (where Girsanov transform appears). It is not clear how to deal
with more general cases as pointed out in \cite[Section 7]{FTHJB}. Moreover
the linearity in the control says that, at the present stage, we cannot
apply this technique to the cases described in Section \ref{APPL1EXIT}.

\item  The regularity assumptions on the coefficients are quite strong and
in any case quite less general than our assumptions. In particular
continuous differentiability is always required also in the Hamiltonian
(assumption 6 above).

\item  The result quoted above is proven for an infinite dimensional case
where mild solutions exist and Girsanov theorem can be applied in a suitable
way. Here we give a finite dimensional version of it. Of course the infinite
dimensional case imposes certain restrictions and probably some technical
assumptions may be removed. However some assumptions that strongly limits
the applicability of the result (linearity in the control, smoothness of the
coefficients) seem to constitute a structural limit of the technique used
there.
\end{itemize}

\textbf{ACKNOWLEDGEMENTS} The authors wish to thank the Referee and the
Editor for the careful reading and for the motivating suggestions.


\end{document}